\documentclass[twocolumn]{autart}    

\usepackage{amsfonts,amsmath,latexsym}
\usepackage{graphicx,color,mathrsfs}
\usepackage[numbers]{natbib}
\usepackage{balance}
\usepackage{hyperref,wasysym}

\newtheorem{proposition}{Proposition} 

\catcode`\@=11
\def\downparenfill{$\m@th\braceld\leaders\vrule\hfill\bracerd$}
\def\overparen#1{\mathop{\vbox{\ialign{##\crcr\crcr \noalign{\kern0.4ex}
\downparenfill\crcr\noalign{\kern0.4ex\nointerlineskip}
$\hfil\displaystyle{#1}\hfil$\crcr}}}\limits}
\catcode`\@=12

\def\power#1#2{\left\lfloor #1\right\rceil^{#2}}

\def\RR{{\mathbb R}}

\def\SS{{\mathbb S}}

\def\gz{\mathfrak{g}}

\def\Ouvs{\mathcal{S}} 

\def\C{\mathcal{C}}
\def\M{\mathcal{M}}
\def\L{\mathcal{L}}
\def\D{\mathcal{D}}

\def\P{\mathcal{P}}

\def\sat{\mathsf{sat}}

\def\diag{\textsf{diag}}
\def\sign{\mathrm{S}}
\def\KR{\mathfrak K} 
\def\w{\mathfrak b} 			
\def\cc{\mathfrak a} 			

\def\vechatZ{\mathbf{\hat z}}
\def\vechatZsol{\mathbf{\hat Z}}
\def\vecZ{\mathbf{z}}
\def\vecZsol{\mathbf{Z}}
\def\vecE{\mathbf{e}}

\def\efin{E}

\def\bare{{\bar e}}
\def\barv{{\bar v}}
\def\bardelta{{\bar \delta }}

\def\hypo{{\sun}} 

\def\bz{\overline{\mbox{\textsf{z}}}}
\def\bu{\overline{\mbox{\textsf{u}}}}
\def\bw{\overline{\mbox{\textsf{w}}}}
\def\bPhi{{\mathchoice%
{\overline{\mbox{\textsf{\char"08}}}}
{\overline{\mbox{\textsf{\char"08}}}}
{\overline{\mbox{\scriptsize\textsf{\char"08}}}}
{\overline{\mbox{\tiny\textsf{\char"08}}}}
}}
\def\bz{{\mathchoice%
{\overline{\mbox{\textsf{z}}}}
{\overline{\mbox{\textsf{z}}}}
{\overline{\mbox{\scriptsize\textsf{z}}}}
{\overline{\mbox{\tiny\textsf{z}}}}
}}

\definecolor{LPcomcolor}{rgb}{1,0,0}
\definecolor{PBcomcolor}{rgb}{1,0,1}
\definecolor{VAcomcolor}{rgb}{0,0,1}
\def\startrevision{\begingroup\color{black}}
\def\stoprevision{\endgroup}

\begin{document}

\begin{frontmatter}

\title{Observers for a non-Lipschitz triangular form} 

\author[First]{P. Bernard}\ead{pauline.bernard@mines-paristech.fr},   
\author[First]{L. Praly}\ead{laurent.praly@mines-paristech.fr},     
\author[Second]{V. Andrieu}\ead{vincent.andrieu@gmail.com}   

\address[First]{Centre Automatique et Syst\`{e}mes, MINES ParisTech, PSL Research University, France }  
\address[Second]{Universit\'e Lyon 1, Villeurbanne, France -- CNRS, UMR 5007, LAGEP, France}


\begin{abstract}                        
We address the problem of designing an observer for triangular non locally Lipschitz dynamical systems. We show the convergence with an arbitrary small error of the classical high gain observer in presence of nonlinearities verifying some H\"older-like condition. Also, for the case when this H\"older condition is not verified, we propose a novel cascaded high gain observer. Under slightly more restrictive assumptions, we prove the convergence of an homogeneous observer and of its cascaded version with the help of an explicit Lyapunov function.

\end{abstract}

\end{frontmatter}

\section{Introduction}

A preliminary step is often 
required in the construction of observers for controlled nonlinear systems. 
It consists in finding a reversible coordinate 
transformation, allowing us to rewrite the system dynamics in a 
target form 
more favorable for writing and/or analyzing the observer.
For example,
the dynamics of a controlled single output system of dimension $n$ which is uniformly observable (see \cite[Definition  I.2.1.2]{GauKup})
and
differentially observable of order $m$ (see \cite[Definition  I.2.4.2]{GauKup}) with $m=n$ can be written with
appropriate coordinates in a Lipschitz triangular form appropriate for
the design of a high gain observer (\cite{GauBor,GauHamOth}).  Such a property is no more true when the order $m$ is strictly larger
than the dimension $n$.  Indeed in this case, we may still get the
usual triangular form but with functions that may not be
Lipschitz (\cite{BerPraAndAUT20161}). 
\startrevision
A particular case of this is when there is only one nonlinear function (in the last line for the single 
output case). This is the so-called phase-variable form.
It has been known for a long time, in particular in the context of dirty-derivatives and output 
differentiation, that a high gain observer can provide an arbitrary small error as long as the nonlinearity is bounded (\cite{Tor} among many others). We also know since \cite{Lev2001} that a sliding mode observer can achieve finite-time convergence under the same assumption. In this paper, we want to build observers for the more general triangular canonical form where non-Lipschitz triangular nonlinearities can appear on any line. As far as we know, this form has not received
much attention apart from its well-known Lipschitz version and the convergence results holding for the phase-variable form do not extend trivially.
\stoprevision

This paper follows and completes \cite{BerPraAndNOLCOS}.  We show here that the classical high gain
observer may still be used when the nonlinearities verify some
H\"older-type condition.  Nevertheless, the asymptotic convergence is lost
and only a convergence with an arbitrary small error remains.  When
the nonlinearities do not verify the required H\"older regularity, it is also
possible to use a cascade of high gain observers, but once again, the
convergence is only with an arbitrary small error.

Fortunately, moving to a generalization of high
gain observer exploiting homogeneity makes it possible to achieve
convergence.
It is at the beginning of the century that researchers started to
consider homogeneous observers with various motivations: exact
differentiators (\cite{Lev2001,Lev2003,Lev2005}), domination as a tool
for designing stabilizing output feedback (\cite{YanLin}, \cite{Qia},
\cite{QiaLin}, \cite{AndPraAst} and references therein (in particular
\cite{AndPraAstCDC2006})), ...  The advantage of this type of
observers is their ability to face H\"older nonlinearities.

\startrevision
With the tools introduced in \cite{AndPraAstSIAM}, we have at our disposal a Lyapunov design to obtain
an homogeneous
observer with degree in $]-1,0[$ for the triangular form
mentioned above. By construction, convergence is guaranteed if the nonlinearities verify a H\"older-type
condition.  We show here that the same Lyapunov design can be extended to the case
where the degree of homogeneity is $-1$. This is interesting since
the constraints on the nonlinearities become less and
less restrictive when the degree gets closer to $-1$.
It turns out that, in the
absence of nonlinearities, the observer we obtain is actually the
exact differentiator presented in \cite{Lev2001} and which
is defined by an homogeneous differential inclusion. But as opposed to \cite{Lev2001} where convergence is 
established via a solution-based analysis, in our case, convergence is again guaranteed by construction since the 
design gives also
a  homogeneous
strict Lyapunov function. Moreover this function enables us to quantify the effect of the observer parameters
on the behavior in presence of H\"older
nonlinearities and disturbances.  Of course, knowing the convergence of the exact differentiator
from \cite{Lev2001}, we could have deduced the existence of such a
Lyapunov function via a converse theorem as in \cite{NakYamNis}. But with only existence, effect 
quantifications as mentioned above is nearly impossible. Actually many efforts have been made to get expressions for 
Lyapunov functions but, as opposed to Lyapunov design, Lyapunov analysis is much harder. As far as we 
know, expressions of Lyapunov functions have been obtained this way only for $m\leq 3$. See \cite{OrtSanMor}.
\stoprevision

Finally, to face the unfortunate situation where the nonlinearities
verify none of the above mentioned H\"older type conditions, we
propose a novel observer made of a cascade of homogeneous observers
whose maximal total dimension is $\frac{m(m+1)}{2}$.  We prove that it
converges without requiring anything on the
nonlinearities (except continuity) in the case where the system trajectories and the input are bounded. 

\startrevision All along our paper, we sometimes use stronger assumptions than necessary in order to simplify
the presentation of our results.  We signal them to the reader with a $\hypo$ symbol
as in \textit{``the trajectories are complete$^\hypo$\  ''}. We discuss how they can be relaxed
later in Section \ref{sec_relaxHypo}, in particular when we restrict our attention to compact sets.
\stoprevision

\textit{Notations}
\begin{enumerate}
	\item
	\label{notation1}
	We define the signed power function as 
	\\[0.7em]\null \hfill $\displaystyle 
	\power{a}{b} = \texttt{sign}(a)\,  |a|^b \, ,
	$\hfill \null \\[0.7em]
	where $b$ is a nonnegative real number. In the particular case where 
	$b=0$, $\power{a}{0}$ is actually any number in the set
	\\[0.7em]\null \hfill $\displaystyle 
	\sign(a)=
	\left\{
	\begin{array}{lc}
	\{1\} & \textrm{if } a>0  \, , \\
	\left[-1,1\right] & \textrm{if } a=0 \, ,\\
	\{-1\} & \textrm{if } a<0 \, .
	\end{array}
	\right.
	$\refstepcounter{equation}\label{LP5}\hfill$(\theequation)$
\\[0.7em]
	Namely, writing $c=\power{a}{0}$ means $c\in\sign(a)$.
	Note that the set valued map $a\mapsto S(a)$ is
	upper semi-continuous with nonempty, compact and 
	convex values.
	\item \label{notation3}
For $(z_1,\ldots,z_i)$ and $(\hat z_{1},\ldots,\hat z_{i})$ (resp. $(\hat z_{i1},\ldots,\hat z_{ii})$) in 
	$\RR^i$, we denote
\\[0.7em]\null \hfill $\displaystyle 
\begin{array}{@{}c@{}}
\vecZ_{{i}}= (z_{1}, \dots, {z}_{i})
\\
\vechatZ_{{i}}=(\hat{z}_{1}, \dots, \hat{z}_{i})
\quad  (\text{resp.} \quad
\vechatZ_{{i}}=(\hat{z}_{i1}, \dots, \hat{z}_{ii}))
\\
e_{ij}=\hat z_{ij}-z_j
\  ,\quad 
e_{j}=\hat z_{j}-z_j
\; ,\  \vecE_{{i}} = \vechatZ_{{i}}-\vecZ_{{i}}
\  .
\end{array}
$
\end{enumerate}

\section{Continuous triangular form}
Consider a nonlinear system of the form
\\[0.7em]\null \hfill $\displaystyle 
	\left\{\begin{array}{rcl}
		\dot{z}_1&=&z_2 + \Phi_1(u,z_1)+w_1 \\
		&\vdots&  \\
		\dot{z}_i&=&z_{i+1} + \Phi_i(u, z_1,\ldots,z_i)+w_i  \\
		&\vdots&  \\
		\dot{z}_m&=&\Phi_m(u,z)+w_m \\
		 y  &=&  z_1+v 
			\end{array}\right. \ , 
$\hfill \null \refstepcounter{equation}\label{eqCanForm}\hfill$(\theequation)$
\\[0.7em]
where $z$ is the state in $\RR^m$, $y$ 
is a measured output in $\RR$, $\Phi$ is a continuous function which is not assumed to be locally Lipschitz.
$w$ can model either a known or an unknown disturbance on the dynamics and $v$ is an unknown disturbance
on the measurement. 
Given locally bounded measurable time functions
$t\mapsto u(t)$ and $t\mapsto w(t)$, we denote $Z(z,t;u,w)$
 a solution of (\ref{eqCanForm}) going through $z$ at time $0$
 which, to simplify the presentation, is assumed to be defined for all $t\geq 0$
(i.e. the trajectories are complete$^\hypo$).  We are interested in estimating $Z(z,t;u,w)$ knowing $y$ and $u$. 

As mentioned in the introduction, this kind of triangular continuous form,
as we call (\ref{eqCanForm}),
appears when we consider systems which are uniformly observable and differentially observable but with an order larger than the system's dimension. An example is given in Section \ref{sec_ex}.

The only existing observer we are aware of able to cope with 
$\Phi$
 no more 
than continuous is the one presented in \cite{BarBouDje}. Its 
dynamics are described by a differential inclusion (see Appendix \ref{app_Barbot})~:
$$
\dot{\hat{z}}\in F(\hat{z},y,u) 
$$
where $(\hat{z},y,u)\mapsto F(\hat{z},y,u) $ is a set valued map. 
In the disturbance free context (i.e. $v=w_i=0$), it 
can be shown that any absolutely continuous solution gives in finite 
time an estimate of $z$ under the only assumption of boundedness of the input and  of the state trajectory. 
But the set valued map $F$ above does 
not satisfy the usual basic assumptions (upper semi-continuous with compact and 
convex values)
(see \cite{Fil,Smi}). It follows that we are not
guaranteed of the 
existence of absolutely continuous solutions nor of possible 
sequential compactness of such solutions and therefore of 
possibilities of approximations of $F$. That is why, in this paper, we look for other candidate observers for the triangular form  \eqref{eqCanForm}.

In doing so, we might have to restrict the possible nonlinearities allowed to obtain the existence of an observer.
The restriction we will impose can be described as follows.
For a positive real number $\cc$,  and  a 
vector  $\alpha$ in $[0,1]^{\frac{m(m+1)}{2}}$, 
we will say that the function $\Phi$ verifies the property  $\P(\alpha, \cc)$ if :

\startrevision
\textbf{Property}   $\boldsymbol{\P(\alpha, \cc)}^\hypo$~:
\bgroup
\itshape
For all $i$ in $\{1,\ldots,m\}$,
for all $\vecZ_{{ia}}$  and $\vecZ_{{ib}}$ in $\RR^m$
and $u$ in 
$U$,  we have\footnote{%
Actually $\Phi_i$ can depend also on $z_{i+1}$ to $z_m$ as long as 
(\ref{eqAss_pertNew}) holds. It can also depend on time requiring some uniform property (see Section \ref{sec_relaxHypo}). 
}~:
\begin{equation}
\left| \Phi_i(u,\vecZ_{{ia}})-\Phi_i(u,\vecZ_{{ib}})\right|\leq \cc \sum_{j=1}^i|z_{ja}-z_{jb}|^{\alpha_{ij}} \ .
\label{eqAss_pertNew}
\end{equation}
\egroup
\stoprevision
This property captures many possible contexts.  In the case in which
 $\alpha_{ij}>0$, it implies that the function $\Phi$ is
H\"older with power $\alpha_{ij}$.  When the $\alpha_{ij}=0$, it
simply implies that the function $\Phi$ 
 is 
bounded.  In the
following,  our aim is to
design an observer depending on the values of $\alpha$.

\startrevision
It is possible to employ the degree of freedom given in (\ref{eqCanForm}) by the time functions
$w$ to deal with the case in which the given function $\Phi(u,z)$ doesn't satisfy $\P(\cc,\alpha)$.  In this
case, an approximation procedure can be carried out to get a function $\hat \Phi$ satisfying $\P(\cc,\alpha)$
and selecting $w = \Phi(u,z)-\hat \Phi(u,z)$ which is an unknown disturbance.  The quality of the estimates
obtained from the observer will then depend on the quality of the approximation (i-e the norm of $w$).  This
is what is done for example in \cite{MorVar} when dealing with locally Lipschitz approximations.
We will further discuss in Section \ref{sec_relaxHypo} how to relax assumption $\P(\cc,\alpha)$.
\stoprevision

In Section \ref{sec_HG}, we start by showing the convergence with an
arbitrary small error of the classical high gain observer when the
nonlinearity $\Phi$ verifies the property $\P$ for certain values of
$\alpha_{ij}$.  We deduce in Section \ref{sec_HG_PR} the convergence
with an arbitrary small error for a cascaded high gain observer when
the input  and the state trajectories are bounded. On an other hand, in
Section \ref{sec_HM}, we show that replacing the high gain structure
by an homogeneous structure enables to obtain convergence under a
slightly more restrictive H\"older restriction.  Then, a cascaded
homogeneous observer is presented in Section \ref{sec_PR}, which
ensures asymptotic convergence when the input  and the state trajectories are bounded. 
As already mentioned, in Section \ref{sec_relaxHypo},
we indicate how the assumptions, marked with $^\hypo$ in the text, can be relaxed.
Finally,
we illustrate our observers with an example in Section \ref{sec_ex}.


\section{High gain observer}
\label{sec_HG}

We consider in this section a classical high gain observer:
\begin{equation}\label{eq_ObsHG}
\left\{\begin{array}{@{}rcl@{}}
{\dot{\hat z}}_{1}&=&{\hat z}_{2}  +  \Phi_1(u,{\hat z}_{1})  
+ \hat{w}_1
- L \, k_{1}\,  (\hat z_{1}-y)
\\
{\dot{\hat z}}_{2}  &=&  {\hat z}_{3} 
+   \Phi_2(u,{\hat z}_{1},{\hat z}_{2})  
+ \hat{w}_2
- L^2 \, k_{2}\,  (\hat z_{1}-y)
\\
& \vdots
\\
{\dot{\hat z}}_{m}&=&
   \Phi_m(u,{\hat z}) 
+ \hat{w}_m
-L^m \, k_{m}\,  (\hat z_{1}-y)
\end{array}\right.
\end{equation} 
where $L$ and the $k_i$'s are gains to be tuned, $y$ is the 
measurement.
The $\hat{w}_i$ are approximations of 
the $w_i$. In particular, when $w_i$ represents unknown 
disturbances, the corresponding $\hat{w}_i$ is simply taken 
 equal to 
0.
In the following, we denote
$$
\Delta w = \hat{w}-w
\  .
$$

When $\Phi$ satisfies the property
$\P(\alpha, \cc)$ with
 $\alpha_{ij}=1$ for all $1\leq j\leq i\leq m$, 
we recognize the usual triangular Lipschitz property
 for which the nominal high-gain observer gives an input to state stability (ISS) property
 with respect to the measurement disturbance $v$ and dynamics disturbance  $w$. 
It is well known that the ISS  gain between the disturbance and the estimation error depends on the high-gain parameter $L$.
Specifically, we have the following
well known result.
See for instance \cite{KhaPra} for a proof.

\startrevision
\begin{proposition}[Nominal high-gain]\label{prop_HG_Nom} 
There exist real numbers $k_1, \dots,  k_m$, $L^*$,  
$\lambda$, $\beta$ and $\gamma$ such that,
\begin{list}{}{%
\parskip 0pt plus 0pt minus 0pt%
\topsep 0pt plus 0pt minus 0pt
\parsep 0pt plus 0pt minus 0pt%
\partopsep 0pt plus 0pt minus 0pt%
\itemsep 0pt plus 0pt minus 0pt
\settowidth{\labelwidth}{a)}%
\setlength{\labelsep}{0.5em}%
\setlength{\leftmargin}{\labelwidth}%
\addtolength{\leftmargin}{\labelsep}%
}
\item[a)]
for all functions $\Phi$ satisfying$^{\hypo}$ for all $i$ and for all $\vecZ_{{ia}}$  and $\vecZ_{{ib}}$ in $\RR^m$
\\[0.5em]\null \hfill $\displaystyle 
\left| \Phi_i(u,\vecZ_{{ia}})-\Phi_i(u,\vecZ_{{ib}})\right| 
\leq \cc \sum_{j=1}^i|z_{ja}-z_{jb}| + \w_i \ 
$\refstepcounter{equation}\label{eq_LipscPerturbed}\hfill$(\theequation)$
\item[b)]
 for all $L\geq \max\{\cc L^*,1\}$,
\item[c)]
for all
locally  bounded time function $(u,v,w,\hat{w})$,
all  $(z,\hat z)$ in
$\RR^m\times\RR^m$, 
\end{list}
any solution $\hat Z(\hat z,z,t;u,v,w,\hat{w})$ 
of (\ref{eq_ObsHG}) verifies, for all $t_0$ and $t$ such that $t\geq t_0\geq 0$, and for all $i$ in $\{1,...,m\}$,
\\[0.7em]\vbox{\hsize=\linewidth\noindent
$\displaystyle 
\left|\hat Z_i(t)-Z_i(t))\right|
$\hfill \null \refstepcounter{equation}\label{eq_ISS_HG}$(\theequation)$\\\null   $\displaystyle 
\leq 
\max\left\{ 
\vrule height 1.2em depth 1.2em width 0pt
L^{i-1} \beta \left|\hat Z_i(t_0)-Z_i(t_0))\right|e^{-\lambda 
L (t-t_0)}, 
\right.$\hfill \null \\[-1em]\null \hfill   $\displaystyle 
\gamma \sup_{\stackrel{1\leq j\leq m}{s\in[t_0,t]}} \left\{  L^{i-1} \, |v(s)| 
, \frac{|\Delta w_j(s)| + \w_j}{L^{j-i+1}}\right\}
\left.\vrule height 1.2em depth 1.2em width 0pt
\hskip -0.5em\right\}
.
$}
where we have used the abbreviation $Z(t)=Z(z,t;u,w)$ and $\hat{Z}(t)=\hat{Z}(z,\hat{z},t;u,v,w,\hat{w})$.
\end{proposition}
\stoprevision

 Since the nominal high-gain observer gives asymptotic convergence for Lipschitz nonlinearities, we may wonder
what type of property is preserved when the nonlinearities are only H\"older. In the following proposition,
we show that the usual  high-gain observer can provide an arbitrary small error on the estimate providing
the H\"older orders $\alpha _{ij}$ satisfy the restrictions given in 
Table \ref{tab_holderPowerPract} or Equation (\ref{eq_alpha}).

\begin{figure}[h]
	\null \hfill $
	\renewcommand{\arraystretch}{1.5}
	\begin{array}[t]{@{}lr|c@{\quad }c@{\quad }c@{\quad }c@{\quad }c@{\quad }c} 
&j&
\scriptstyle 1 & \scriptstyle 2 &\ldots & \scriptstyle m-2 & \scriptstyle  m-1 & \scriptstyle m
\\[-0.5em]
i &&
	\\[-0.3em]\hline
	\scriptstyle 1 &
	&\textstyle \frac{m-2}{m-1} &  & &  & &
	\\
\scriptstyle 2&
	&\textstyle\frac{m-3}{m-2} &\textstyle\frac{m-3}{m-2} &  &  & &
	\\
\vdots&\quad  \alpha _{ij} > \  
	&\vdots&\vdots&\ddots&&&
	\\
\scriptstyle m-2&
	&\textstyle\frac{1}{2}&\textstyle\frac{1}{2}&\ldots&\textstyle\frac{1}{2}&&
	\\
\scriptstyle m-1&
	&\textstyle0&\textstyle0&\ldots&\ldots&\textstyle0&
	\\\hline
\scriptstyle m&\quad  \alpha _{mj} \geq   \  
	&0&0&\ldots&\ldots&\ldots&0
	\end{array}
	$\hfill \null\\[0.5em]\null \hfill  
	Table \refstepcounter{table}\arabic{table} 
	\label{tab_holderPowerPract}: H\"older restrictions on $\Phi$ 
	for arbitrarily small errors with a high gain observer.\hfill \null 
\end{figure}

\begin{proposition}\label{prop_HG}
Assume the function $\Phi$  verifies $\P(\alpha,\cc)$ for some 
$(\alpha,\cc)$ in $[0,1]^{\frac{m(m+1)}{2}}\times \RR_+$ satisfying, for $1\leq j\leq i$
\begin{equation}
\label{eq_alpha}
\begin{array}{r@{\quad }c@{\quad }l} 
\frac{m-i-1}{m-i}  < \alpha_{ij} \leq 1 & \text{for} & i=1\dots,m-1\ ,\\
0\leq \alpha_{mj}\leq 1
\end{array}
\end{equation}
Then, there exist real numbers $k_1, \dots,  k_m$, such that,  for all 
$\epsilon>0$ we can find positive real numbers $\lambda$, $\beta$, 
$\gamma$, and $L^*$ such that, 
for all $L\geq L^*$,  for all
locally  bounded time function $(u,v,w,\hat{w})$
 and all  $(z,\hat z)$ in
$\RR^m\times\RR^m$,  
any solution $\hat Z(\hat z,z,t;u,v,w,\hat{w})$ 
of (\ref{eq_ObsHG}) verifies, for all $t_0$ and $t$ such that $t\geq t_0\geq 0$, and for all $i$ in $\{1,...,m\}$,
\\[0.7em]\vbox{\hsize=\linewidth\noindent$\displaystyle 
\left|\hat Z_i(t)-Z_i(t))\right|
$\hfill \null \\\null   $\displaystyle 
\leq 
\max\left\{ 
\vrule height 1.2em depth 1.2em width 0pt
\right.
\epsilon \,  ,\,  
L^{i-1} \beta \left|\hat Z_i(t_0)-Z_i(t_0))\right|e^{-\lambda L (t-t_0)}, 
$  \null \\\null \hfill   $\displaystyle 
\gamma \sup_{\stackrel{1\leq j\leq m}{s\in[t_0,t]}} \left\{  L^{i-1} \, |v(s)| 
, \frac{|\Delta w_j(s)|}{L^{j-i+1}}\right\}
\left.\vrule height 1.2em depth 1.2em width 0pt\right\}
$}
where we have used the abbreviation $Z(t)=Z(z,t;u,w)$ and $\hat{Z}(t)=\hat{Z}(z,\hat{z},t;u,v,w,\hat{w})$.
\end{proposition}
\par\vspace{1em}\noindent
Comparing this inequality with (\ref{eq_ISS_HG}), we have now 
the arbitrarily small non zero $\varepsilon $ in the right hand 
side but this is obtained under the H\"{o}lder 
condition instead of the Lipschitz one.

\startrevision
\begin{pf}
With Young's inequality, we obtain from (\ref{eqAss_pertNew}) that, for all $\sigma_{ij}$ in $\RR_+$
and all $\vechatZ$ and $\vecZ$ in $\RR^m$
\\[0.5em]\null \hfill $\displaystyle 
\left| \Phi_i(u,\vechatZ_{{i}})-\Phi_i(u,\vecZ_{{i}})\right|
\; \leq  \; \sum_{j=1}^i\cc_{ij} |\hat{z}_j-z_j| + \w_{ij}
\  ,
$\refstepcounter{equation}\label{eq_Young}\hfill$(\theequation)$\\[0.5em]
with $\cc_{ij}$ and $\w_{ij}$ defined as
\begin{equation}
\label{LP6}
\!\left\{ \begin{array}{@{}l@{\quad }l@{}}
\cc_{ij} = 0 \ ,\ \w_{ij}=\cc\ ,\ & \text{ if } \alpha_{ij}=0\\
\cc_{ij} = \cc^{\frac{1}{\alpha_{ij}}} \alpha_{ij}\sigma_{ij}^{\frac{1}{\alpha_{ij}}}\ ,\ \w_{ij} =   \frac{1-\alpha_{ij}} { \sigma_{ij}^{
\frac{1}{1-\alpha_{ij}}
}
} 
& \text{ if }0<\alpha_{ij}<1\\
\cc_{ij} = \cc\ ,\ \w_{ij} =0&\text{ if }\alpha_{ij}=1
\end{array}
\right.
\end{equation}
With (\ref{eq_Young}), the assumptions of Proposition 
\ref{prop_HG_Nom} are satisfied with $\w_i =  \sum_{j=1}^i\w_{ij}$.
It gives
$k_1, \dots,  k_m$, $L^*$, 
$\lambda$, $\beta$ and $\gamma$ and,
if $L>\max_{i\geq j} \left\{\cc_{ij}L^*,1 \right\}$,
the solution satisfies
the ISS inequality (\ref{eq_ISS_HG}).
The result will follow if there exist $L$ and $\sigma_{ij}$ such that 
\\[0.5em]\null \hfill $\displaystyle 
L>\max_{i\geq j} \left\{\cc_{ij}
 L^*,1\right\}
\    ,\quad  
\max_{i, j} \sum_{\ell=1}^j
\gamma\w_{j\ell}
L^{i-j-1}
\leq \epsilon\ .
$\refstepcounter{equation}\label{eq_ConstrL}\hfill$(\theequation)$ \\[0.5em]
At this point, we have to work with the expressions of $\cc _{ij}$ and $\w_{j\ell}$ given in (\ref{LP6}).
From (\ref{eq_alpha}), $\alpha _{ij}$ can be zero only if $i=m$. And, when $\alpha _{m\ell}=0$, we get
\\[0.5em]\null \hfill $\displaystyle 
\gamma \w_{m\ell}L^{i-m-1}
\;=\; \gamma \cc L^{i-m-1}\; \leq \; \frac{\gamma \cc}{L}
$\hfill \null \\[0.5em]
Say that we pick $\sigma _{m\ell}=1$ in this case.
For all the other cases, we choose
\\[0.5em]\null \hfill $\displaystyle 
\sigma_{j\ell}=\left(\frac{2j\gamma}{\epsilon}  (1-\alpha_{j\ell})L^{(m-j-1)}\right)^{1-\alpha_{j\ell}} 
\  ,
$\hfill \null \\[0.5em]
to obtain from (\ref{LP6})
\\[0.5em]\null \hfill $\displaystyle 
\gamma \w_{j\ell}L^{i-j-1}\; \leq \; \epsilon\frac{1}{j}\frac{ 1}{2L^{m-i}}
$\hfill \null \\[0.5em]
So, with this selection of the $\sigma _{j\ell}$, the right inequality in (\ref{eq_ConstrL}) is satisfied for $L$ sufficiently large.
Then, according to (\ref{LP6}), the $a_{ij}$ are independent of $L$ or proportional to 
$L^{(m-i-1)\frac{1-\alpha_{ij}} {\alpha _{ij}}}$. But 
with \eqref{eq_alpha} we have 
\\[0.5em]\null \hfill $\displaystyle 
0< (m-i-1)\frac{1-\alpha _{ij}}{\alpha _{ij}}
<1 \ .
$\hfill \null \\[0.5em]
This implies that $\frac{\cc _{ij}}{L}$ tends to $0$ as $L$ tends to $+\infty $.
We conclude that (\ref{eq_ConstrL}) holds if we pick $L$ sufficiently large.
\end{pf}
\stoprevision

It is interesting to remark the weakness of the assumptions imposed on the last two components of the function
$\Phi$.  Indeed, (\ref{eq_alpha}) only imposes that $\Phi_{m-1}$ be H\"older without any restriction on the
order, and that $\Phi_m$ be bounded$^\hypo$.

\section{Cascaded high gain observer}
\label{sec_HG_PR}

According to Proposition \ref{prop_HG}, the classical high gain observer can provide an arbitrary small error 
when the last nonlinearity is only  bounded and when there is no disturbance.
We exploit here this observation by proposing the following cascaded high gain observer
to deal with the case where the functions $\Phi_i$ do not satisfy 
(\ref{eq_alpha}):
\vbox{$$
\begin{array}{@{}rcl@{}}
{\dot{\hat z}}_{11}&=& \hat{w}_1 - L_1 \, k_{11}\,  (\hat z_{11}-z_1)
\\[-0.7em]
\multicolumn{3}{@{}c@{}}{\dotfill}\\
{\dot{\hat z}}_{21}&=&{\hat z}_{22} + \Phi_1 (u,{\hat z}_{11}) 
+ \hat{w}_1
- L_2 \, k_{21}\,  (\hat z_{21}-z_1)
\\
{\dot{\hat z}}_{22}&=& \hat w_2 - L_2^2 \, k_{22}\,  (\hat z_{21}-z_1)
\\[-0.7em]
\multicolumn{3}{@{}c@{}}{\dotfill}\\
&
\vdots
&
\\[-0.7em]
\multicolumn{3}{@{}c@{}}{\dotfill}\\
{\dot{\hat z}}_{m1}&=&{\hat z}_{m2}  + \Phi_1  (u,{\hat z}_{(m-1)1})
+ \hat{w}_1 
- L_m \, k_{m1}\, (\hat z_{m1}-z_1)
\\
{\dot{\hat z}}_{m2}&=&{\hat z}_{m3}  +\Phi_2  (u,{\hat z}_{(m-1)1},{\hat 
	z}_{(m-1)2}) \\
&&\hfill+ \hat{w}_2
- L_m^2 \, k_{m2}\,  (\hat z_{m1}-z_1)
\\
& \vdots
\\
{\dot{\hat z}}_{mm}&=& \hat w_m
- L_m^m \, k_{mm}\,  (\hat z_{m1}-z_1)
\end{array}
$$
\\[-2.6em]\null \hfill 
\refstepcounter{equation}\label{eqObsHG_PR}$(\theequation)$}\\[1em]
with the gain $k_{ij}$ chosen as in a classical high gain observer of dimension $i$,
$\hat{w}_i$ are estimations of $w_i$ and $L_i$ are the high gains parameters to be chosen. 

Assuming the input function and the system solution are bounded, it is shown in the following that estimation
with an arbitrary small error can be achieved by the cascaded high-gain observer.
\startrevision
\begin{proposition}
\label{prop_HG_PR}
Assume $\Phi$ is continuous. 
For any positive real numbers $\bz$ and $\bu$, for any strictly positive real number $\epsilon$, there exist a choice of $(L_1,...,L_m)$, a class $\mathcal {KL}$
function $\beta$ and two class $\mathcal K_\infty$ functions 
$\gamma_1$ and $\gamma_2$ such that,
 for all 
locally  bounded time function $(u,v,w,\hat{w})$,
for all  $(z,\hat z)$ in
$\RR^m\times\RR^m$ and  for all $t$ such that $|Z(z,s;u,w)|\leq \bz$ and $|u(s)|\leq \bu$ for all $0\leq s\leq t$, 
any solution $\left(\vechatZsol_1(\hat z,z,t;u,v,w,\hat{w}),...,\vechatZsol_m(\hat z,z,t;u,v,w,\hat{w})\right)$
 of (\ref{eqObsHG_PR}) verifies, for all $i$ in $\{1,\ldots,m\}$,
\\[0.6em]\vbox{\hsize=\linewidth\noindent$\displaystyle 
|\mathbf{\hat Z} _{{i}}(t)-\mathbf{Z}_i(t)|
$\hfill \null \\[-0.5em]\null \quad  $\displaystyle 
\leq 
\max\left\{ 
\vrule height 1.2em depth 1.2em width 0pt
\right.\hskip -0.3em
\varepsilon \,  ,\,  
\beta\left(\sum_{j=1}^i|\hat{z}_j-z_j| ,t\right) , 
$\hfill \null \\[-0.7em]\null \hfill   $\displaystyle 
\sup_{s\in [0,t]}\Big\{\gamma_1(|v(s)|), \gamma_2(|\Delta w(s)|)\Big\}
\hskip -0.5em
\left.\vrule height 1.2em depth 1.2em width 0pt\right\}
$}
where $\vechatZsol_i$ is the state of the $i$th block (see Notation \ref{notation3}) and we have used the abbreviation
$\vechatZsol_i(t)=\vechatZsol_i(\hat z,z,t;u,v,w,\hat{w})$ and $\vecZsol_i(t)=\vecZsol_i(z,t;u,w)$.
\end{proposition}
\stoprevision

\begin{pf}
This result is nothing but a straightforward consequence of the fact 
that a cascade of ISS systems is ISS.

Specifically the error system 
attached to the high gain observer in block $i$ has state 
$\vecE_{{i}}$ (see Notation \ref{notation3}) and input
$v$ and
$\delta _{ij}$ defined as
$$
\begin{array}{rcl}
\delta _{ij}&=&\left[
\Phi_j (u,\vechatZ_{{(i-1)}})-\Phi_j(u,\vecZ_{{(i-1)}})\right]+\left[\hat{w}_j-w_j
\right]
\\
\delta _{ii}&=& -z_{i+1}-\Phi_i(u,\vecZ_{{i}})+\hat{w}_i-w_i
\end{array}
$$
with $z_{m+1}=0$.
With  Proposition \ref{prop_HG_Nom}, we have the existence of
$k_{i1}, \dots,  k_{ii}$, 
$\lambda _i$, $\beta _i$ and $\gamma _i$ such that we have, for all 
$L_i\geq 1$,  all $t\geq t_i\geq 0$, all $j$ in $\{1,\ldots,i\}$ and with $e_{ij}(t)$ denoting the $j$th error in the $i$th
block evaluated along the solution at time $t$,
\\[0.7em]\vbox{\hsize=\linewidth\noindent$\displaystyle 
\left|e_{ij}(t)\right|
$\hfill \null \\\null \quad  $\displaystyle 
\leq 
\max\left\{ 
\vrule height 1.2em depth 1.2em width 0pt
L_i^{j-1} \beta _i\left|\vecE_{{i}}(t_i)\right|
e^{-\lambda _i L_i (t-t_i)}, 
\right.
$\hfill \null \\[-1em]\null \hfill   $\displaystyle 
\gamma _i\sup_{\stackrel{1\leq \ell\leq j}{s\in[t_i,t]}}
\left\{  L_i^{j-1} \, |v(s)| 
, \frac{|\delta _{i\ell}(s)|}{L_i^{\ell-j+1}}
\right\}
\left.\vrule height 1.2em depth 1.2em width 0pt\right\}
\  .
$}\\[0.7em]
But  the  continuity   of the $\Phi_j$ 
 implies the existence of a function\footnote{%
Simply take $\rho(s) = \max_{|u|\leq \bu, |\vecZ_j|\leq \bz,  |e|\leq 
s}|\Phi_j(u,\vecZ_j+e)-\Phi_j(u,\vecZ_j)|$.}
${\rho}$ of class $\mathcal{K}$ such that, for all $j$ in $\{1,\dots,m\}$ and for all
$(\vecZ_{{(i-1)}},\vechatZ_{{(i-1)}},u)$ 
in  $\RR^{i-1}\times\RR^{i-1}\times U$  satisfying $|\vecZ_{{(i-1)}}|\leq \bz$ and $|u|\leq \bu$,
$$
| \Phi_j(u,\vechatZ_{{(i-1)}})-\Phi_j(u,\vecZ_{{(i-1)}})|\leq 
\rho\left(|\vecE_{{(i-1)}}|\right) \ .
$$
This implies
\begin{eqnarray*}
|\delta _{i\ell}(s)| &\leq &\rho(|\vecE_{{i-1}}(s)|) + |\Delta w_\ell(s)|\ ,\ \ell=1,\dots,j-1\ ,\\
|\delta _{ii}(s)| &\leq &\bz_{i+1}+\bPhi _i + |\Delta w_i(s)| \ ,
\end{eqnarray*}
where $\bPhi _i=\max_{|u|\leq \bu,|\vecZ_i|\leq \bz} \left|\Phi_i(u,\vecZ_i)\right|$.
Hence, we have the existence of $c_i$ independent of $L_i$ 
	such that
\\[0.7em]\vbox{\hsize=\linewidth\noindent$\displaystyle 
|\vecE_{{i}}(t)| 
$\hfill \null \\\null   $\displaystyle 
\leq c_i
\max\left\{ 
\vrule height 1.2em depth 1.2em width 0pt
L_i^{i-1}  \left|\vecE_{{i}}(t_i)\right|
e^{-\lambda _i L_i (t-t_i)}
\,  ,\,  
\sup_{s\in[t_i,t]}
L_i^{i-1} \, |v(s)| 
\,  ,
\right.
$\hfill \null \\\null \hfill   $\displaystyle 
\sup_{s\in[t_i,t]}
\frac{ \rho(|\vecE_{{i-1}}(s)|)}{L_i^{2-i}}
\,  ,\,  
\sup_{
\renewcommand{\arraystretch}{0.8}
\begin{array}{@{}c@{}}
\scriptstyle
1\leq \ell\leq i
\\
\scriptstyle
s\in[t_i,t]
\end{array}}
 \frac{|\Delta w_\ell(s)|}{L_i^{\ell-i+1}}
\,  ,\,  
\frac{\bz_{i+1}+\bPhi _i }{L_i}
\left.\vrule height 1.2em depth 1.2em width 0pt\right\}
\  .
$}\\[0.7em]
This makes precise what we wrote above that we have a cascade of ISS systems.
Hence  (see \cite[Prop. 7.2]{Son}), for each $i$ in $\{1,\dots,m\}$, there exist a class 
$\mathcal{KL}$ function $\bar \beta_i$ and class $\mathcal K$ 
functions $\gamma_{vi}$ and $\gamma_{wi}$, each depending on $L_1$ 
to $L_i$ and such that we have, for all $t\geq 0$,
\\[0.7em]
\vbox{\hsize=\linewidth\noindent$\displaystyle 
|\vecE_{{i}}(t)| \leq
\max\left\{ 
\vrule height 1.2em depth 1.2em width 0pt
\bar\beta_i\left(\max_{j\in\{1,\ldots,i\}}\{|\vecE_{{j}}(0)|\},t\right)\,  ,
\right.
$\hfill \null \\\null \hfill $\displaystyle 
\left.
\varpi _i\,  ,\,  
\sup_{s\in[0,t]}\left\{\gamma_{vi}( |v(s)| ) , \gamma_{wi}( |\Delta w(s)|)\right\}
\vrule height 1.2em depth 1.2em width 0pt
\right\}\ .
$}\\[0.7em]
where $\varpi _i$ is a positive real number defined by the sequences
$$
\varpi _1 = c_1\frac{\bz_{2}+\bPhi  _1}{L_1}\ ,\ \varpi _{i} = c_i\max\left\{\frac{\bz_{i+1}+\bPhi _i }{L_i}, \frac{\rho(\varpi _{i-1})}{L_i^{2-i}}\right\} \ .
$$
Then by picking $L_i\geq L_i^*$ where $L_i^*$ is defined recursively as~:
$$
\begin{array}{c}
    \epsilon_m = \epsilon \ , \quad\epsilon_i = 
    \min\left(\epsilon, 
    \rho^{-1}\left(\frac{\epsilon_{i+1}}{c_{i+1}L_{i+1}^{i-2}}\epsilon_{i+1}\right)\right)
\\\displaystyle 
 L_m^*= \frac{c_m\bPhi _m}{\varepsilon _m}
\  ,\quad 
L_i^*= \frac{c_i[\bz_{i+1}+\bPhi  _i]}{\varepsilon _i}
\end{array}
$$
we obtain $\varpi _i\leq \epsilon$ for all $i$, hence the result.
\end{pf}

This observer has the advantage of working without any assumption on 
the nonlinearities besides their continuity. Note however that it requires the knowledge of a bound on the 
system solution and on the input.
  Also we
may not need to build $m$ blocks, since according to
Proposition \ref{prop_HG}, we need to create a new block only for the
indexes $i$ where $\Phi_i$ does not verify Property  $\P(\alpha,\cc)$
for any $\cc\geq 0$ and with $\alpha$ satisfying (\ref{eq_alpha}).
Unfortunately, as it appears from the proof of Proposition \ref{prop_HG_PR}, the
choice of $(L_1,...,L_m)$ can be complicated.
Besides, only a convergence with an arbitrary small error is 
obtained. It may thus be necessary to take very high gains which is 
problematic in terms of peaking and most importantly in presence of 
noise (see Section \ref{sec_ex}). In the following two sections, we 
move our attention to homogeneous observers, and show that they 
enable to obtain convergence.


\section{Homogeneous observer}
\label{sec_HM}

Homogeneous observers are extensions of high gain observers able to
cope with some non Lipschitz functions.  As mentioned in the
introduction, they already have an old history (see \cite{Lev2001},
\cite{Lev2003}, \cite{YanLin}, \cite{Lev2005}, \cite{Qia},
\cite{QiaLin}, \cite{AndPraAstCDC2006}, \cite{AndPraAstSIAM},
\cite{AndPraAst} ).  In our context they take the form~:
\begin{equation}\label{eq_ObsHomo}
\begin{array}{@{}rcl@{}}
{\dot{\hat z}}_{1}&=&{\hat z}_{2}  +  \Phi_1  (u,{\hat z}_{1},t) 
+\hat{w}_1
-L \, k_{1}\,  \power{\hat z_{1}-y}{\frac{r_2}{r_1}}
\\
{\dot{\hat z}}_{2}  &=&  {\hat z}_{3} 
+  \Phi_2   (u,{\hat z}_{1},{\hat z}_{2},t) 
+\hat{w}_2
-L^2 \, k_{2}\,  \power{\hat z_{1}-y}{\frac{r_3}{r_1}}
\\
& \vdots
\\
{\dot{\hat z}}_{m}&=&
 \Phi_m   (u,{\hat z},t) 
+\hat{w}_m
-L^m \, k_{m}\,  \power{\hat z_{1}-y}{\frac{r_{m+1}}{r_1}}
\end{array}
\end{equation}
\startrevision
where $r$ is a vector in $\RR^{m+1}$, called weight vector, the components of which, called weights, are defined by
\begin{equation}
\label{eqdefr}
r_i = 1-d_0(m-i) \ ,
\end{equation}
\stoprevision
and where $L$ and the $k_i$'s are gains to be tuned, $d_0$ is a parameter to 
be chosen in $[-1,0]$.
We refer to Notation 
(\ref{notation1}) for the case $d_0=-1$,
for which the dynamics \eqref{eq_ObsHomo} must
be understood as a differential inclusion.
When $d_0=0$, we recover the high-gain observer studied in Section \ref{sec_HG}.
As mentioned in Proposition \ref{prop_HG}, the usual high-gain
observer can provide an estimation with an arbitrary small error provided the
nonlinearity satisfies the property $\P(\alpha,\cc)$ with the
$\alpha_{ij}$ verifying (\ref{eq_alpha}).  In the
following proposition we claim that asymptotic estimation may be obtained
with homogeneous correction terms and when considering
nonlinearities which satisfies $\P(\alpha,\cc)$ with the
$\alpha_{ij}$ verifying
\begin{equation}\label{eq_defalphaHomo}
\alpha_{ij}  = 
\frac{1-d_0(m-i-1)}{1-d_0(m-j)}=\frac{r_{i+1}}{r_j}\ ,\ 1\leq j\leq i \leq m\ .
\end{equation}
Those conditions in the extreme case where $d_0=-1$ are summed up in Table \ref{tab_holderPower}. 
On top of that, finite time estimation may be obtained.
\begin{figure}
	\null \hfill $
	\renewcommand{\arraystretch}{1.5}
	\begin{array}[t]{@{}lr|c@{\quad }c@{\quad }c@{\quad }c@{\quad }c@{\quad }c} 
&j&
\scriptstyle 1 & \scriptstyle 2 &\ldots & \scriptstyle m-2 & \scriptstyle  m-1 & \scriptstyle m
\\[-0.5em]
i &&
\\[-0.3em]\hline
\scriptstyle 1 &
&\textstyle \frac{m-1}{m} &  & &  & &
\\
\scriptstyle 2&
&\textstyle\frac{m-2}{m}&\textstyle\frac{m-2}{m-1}&  &  & &
\\
\vdots&\quad  \alpha _{ij} = \  
&\vdots&\vdots&\ddots&&&
\\
\scriptstyle m-2&
&\textstyle\frac{2}{m}&\textstyle\frac{2}{m-1}&\ldots&\textstyle\frac{2}{3}&&
\\
\scriptstyle m-1&
&\textstyle\frac{1}{m}&\textstyle\frac{1}{m-1}&\ldots&\ldots&\textstyle\frac{1}{2}&
\\
\scriptstyle m&
&0&0&\ldots&\ldots&\ldots&0
\end{array}
$\hfill \null \\[0.5em]\null \hfill  
Table \refstepcounter{table}\arabic{table} \label{tab_holderPower}:
H\"older restrictions on $\Phi$ for a homogeneous observer
with $d_0=-1$\hfill \null\\[-0.5em]\null \hrulefill\null 
\end{figure}


\begin{proposition}\label{prop_HM_pert}
Assume that there exist $d_0$ in $[-1,0]$ and  $\cc$ in $\RR_+$ such that
$\Phi$   satisfies  $\P(\alpha,\cc)$ 
with $\alpha$ verifying (\ref{eq_defalphaHomo})$^\hypo$.
There exist $(k_1, \dots,  k_m)$, such that for all $\bar w_m>0$
there exist $L^*\geq 1$ and a positive constant $\gamma$ such that,
for all $L\geq L^*$ there exists a class $\mathcal {KL}$ function
$\beta$ such that
 for all
locally  bounded time function $(u,v,w,\hat{w})$,
and all  $(z,\hat z)$ in
$\RR^m\times\RR^m$ 
 system (\ref{eq_ObsHomo})
admits absolutely continuous solutions $\hat Z(\hat z,z, t;u,v,w,\hat{w})$ defined on
$\RR_+$
and for any such solution the following implications hold for all $t_0$ and $t$ such that
$t\geq t_0\geq 0$, and for all $i$ in $\{1,...,m\}$ :
\\\noindent \underline{If $d_0>-1$~:}
\\[0.7em]$\displaystyle 
	|\hat Z_i(t)-Z_i(t)| \leq 
\max\left\{ 
\vrule height 1.2em depth 1.2em width 0pt
\right.
\beta(|\hat Z(t_0)- Z(t_0)|,t-t_0)\,  ,\,  $\hfill \null 
\refstepcounter{equation}\label{LP1}\hfill$(\theequation)$
\\\null \hfill $\displaystyle 
\gamma \sup_{\stackrel{1\leq j\leq i}{s\in[t_0,t]}} \left\{ L^{i-1}|v(s)|^{\frac{r_i}{r_1}},
\frac{|\Delta w_j(s)|^{\frac{r_i}{r_{j+1}}}}{L^{\mu_{ij}}}\,  \right\} 
\left.
\vrule height 1.2em depth 1.2em width 0pt\right\}$\\[0.7em]
	where $\mu_{ij} = (j-i+1) \frac{r_1}{r_{j+1}}$, and we have used the abbreviation
$Z(t)=Z(z,t;u,w)$ and $\hat{Z}(t)=\hat{Z}(z,\hat{z},t;u,v,w,\hat{w})$.
\\
Moreover, when $d_0<0$ and  $v(t)=w_j(t)=0$ for all $t$ and $j=1,\dots,m$,  
	 there exists $T$	such that
$\hat Z(\hat z,z,t)=Z(z,t)$ for all $t\geq T$.

\noindent \underline{If $d_0=-1$ and  $|\Delta w_m(t)|\leq \bar w_m$~:}
\\[0.7em]$\displaystyle 
	|\hat Z_i(t)-Z_i(t)| \leq 
\max\left\{ 
\vrule height 1.2em depth 1.2em width 0pt
\right.
\beta(|\hat Z(t_0)- Z(t_0)|,t-t_0)\,  ,
$\hfill \null 
\refstepcounter{equation}\label{LP2}\hfill$(\theequation)$
\\\null \hfill $\displaystyle 
\gamma \sup_{\stackrel{1\leq j\leq i-1}{s\in[t_0,t]}} \left\{ L^{i-1}|v(s)|^{\frac{r_i}{r_1}},
 \frac{|\Delta w_j(s)|^{\frac{r_i}{r_{j+1}}}}{L^{\mu_{ij}}}  \,  \right\}
\left.
\vrule height 1.2em depth 1.2em width 0pt\right\}$\\[0.7em]
where $\mu_{ij}$, $Z(t)$ and $\hat{Z}(t)$ are defined above.\\
Moreover, when  $v(t)=w_j(t)=0$  for all $t$ and $j=1,\dots,m$, there exists $T$ such that
$\hat Z(t)=Z(t)$ for all $t\geq T$.
\end{proposition}

Note that $j$ is in $\{1,\ldots,i\}$ in (\ref{LP1}) whereas it is 
in $\{1,\ldots,i-1\}$ in (\ref{LP2}).

The proof of Proposition \ref{prop_HM_pert} for the case $d_0\in]-1,0]$ and without disturbances is given 
for example in \cite{AndPraAstSIAM}.
Actually \cite{AndPraAstSIAM} gives a Lyapunov design of the observer (\ref{eq_ObsHomo}) with a recursive construction of both 
Lyapunov function and observer. 
Here we are concerned with the case $d_0=-1$.
In this limit case, the observer \eqref{eq_ObsHomo} is a differential inclusion corresponding to the exact differentiator
studied in \cite{Lev2001}, where convergence is established in the particular case in which $\Phi_i=0$ for
$j=1,\dots, m-1$ and  $\Phi_m$ is bounded. 
We prove in  Proposition \ref{Lemm_StricteLyap} that the Lyapunov design of \cite{AndPraAstSIAM} can be 
extended to this case. This allows us to
show that the observer \eqref{eq_ObsHomo} still converges if, for each $i$, $\Phi_i$ is  H\"older 
with order 
$\alpha _{ij}$  equal to the values given in Table \ref{tab_holderPower}, where
$i$ is the index of $\Phi_i$ and $j$ is the index of $e_j$. 
We also recover the same bound in presence of a noise $v$ as the one given in \cite{Lev2001}. 

\startrevision
Actually some effort has been devoted to Lyapunov analysis for establishing the convergence of the observer 
proposed in \cite{Lev2001}. But, as far as we are aware of, this more difficult route has been successful for 
$m\leq 3$ only. See \cite{OrtSanMor}.
\stoprevision

Finally, it is interesting to remark that in the case $d_0=-1$ the ISS property between the disturbance $w_m$ and the estimation error
is with restrictions as defined in \cite[Definition 3.1]{Tee}.
If $|\Delta w_m(t)|\leq \bar w_m$ and $L$ is chosen sufficiently large, then asymptotic convergence is obtained. However, nothing can be said when $|\Delta w_m|>\bar w_m$.
Moreover, it may be possible for a bounded large disturbance to induce a norm of the estimation error which goes to infinity. 
We believe that this problem could be solved employing homogeneous in the bi-limit observer as in \cite{AndPraAstSIAM}.
\startrevision
It is shown to be doable in dimension 2 in \cite{CruMorFri}.
\stoprevision

\begin{pf}
The set-valued function
$e_1\mapsto \power{e_1}{0}$ defined in Notations \ref{notation1} is
upper semi-continuous and has convex and compact values.  Thus,
according to \cite{Fil}, there exist absolutely continuous solutions
to \eqref{eq_ObsHomo}.

Let $\L=\diag(1,L,...,L^{m-1})$.
The error $e=\hat{z}-z$ produced by the observer \eqref{eq_ObsHomo} satisfies
\begin{equation}\label{eq_Syste}
\dot e \in LS_m e + \delta + L\L \KR(e_1+v)
\end{equation}
where $S_m$ is the shifting matrix of order $m$,
$$
\delta =  \Phi(u,\hat z)  +\hat w - \Phi(u,z) - w\ ,
$$
and $\KR$ is the homogeneous correction term the components of which are defined as
$$
(\KR(e_1))_i =  k_i\power{e_1}{\frac{r_{i+1}}{r_1}}
$$
where $(k_1, \dots, k_m)$ are positive real number  and $r_i$ is defined in \eqref{eqdefr}.

Let also $V:\RR^m\to \RR_+$ be the function defined as
\begin{multline}\label{eq_defV}
V(\bare) = \sum_{i=1}^{m-1}\int_{\power{\bare_{i+1}}{\frac{r_{i}}{r_{i+1}}}}^{\ell_i \bare_i}
\left[\power{x}{\frac{d_V-r_i}{r_i}}- 
\power{\bare_{i+1}}{\frac{d_V-r_i}{r_{i+1}}} \right]dx \\+\frac{|\bare_m|^{d_V}}{d_V},
\end{multline}
where $d_V$ and $\ell_i$ are positive real numbers such that
$d_V>2m-1$.  Note that $V$ is a homogeneous function with weight vector $r$.
It is nothing but the one proposed in \cite[Theorem 3.1]{AndPraAstSIAM}  for designing an observer 
homogeneous in the bi-limit with $d_0$ in $]-1,0]$. There it is shown that,
by appropriately
selecting the parameters $\ell_i$ and $k_i$, $V$ is a strict $C^1$
Lyapunov function homogeneous of degree $d_V$ 
for the $L$-independent auxiliary 
system with state $\bare$ ~:
\begin{equation}
\label{LP3}
\dot \bare \in S_m \bare  + \KR(\bare_1)\ .
\end{equation}

With this result in hand a robustness analysis can be carried out on a system of the form (\ref{eq_Syste}).
In fact, the same approach can be followed for the case $d_0=-1$ and 
the following technical result is proved in Appendix \ref{Sec_ProofLem1}.
\begin{lem}\label{Lemm_Rob}
For all $d_0$ in $[-1,0]$, the function $V$ defined in (\ref{eq_defV}) is positive definite and there exist positive real numbers $k_1, \dots k_m$, $\ell_1, \dots\ell_m$, $\lambda$, 
$c_\delta$ and $c_v$ such that for all $\bare$ in $\RR^m$, $\bardelta$ in $\RR^m$ and $\barv$ in $\RR$
the following implication  holds~:
\\[0.7em]
if \hfill   $
\displaystyle |\bardelta_i| \leq \displaystyle c_\delta 
V(\bare)^\frac{r_{i+1}}{d_V}
\!,\quad \forall i\,  ,$\hfill and \hfill $
\displaystyle |\barv|\leq \displaystyle c_v V(\bare)^\frac{r_1}{d_V}
$\hfill then\footnote{%
Here the $\max$ is with respect to $s$ in $\power{\bare_1+\barv)}{0}$ 
appearing in the $m$th component $\KR(\bare_1+\barv)_m$ of $\KR(\bare_1+\barv)$.
}\\[0.5em]$\displaystyle 
\max\!\left\{\frac{\partial V}{\partial \bare}(\bare) (S_m(\bare) + \bardelta + \KR(\bare_1+\barv))\right\} 
\leq -\lambda V(\bare)^\frac{d_V+d_0}{d_V}
\hskip-0.5em .
$
\end{lem}

This Lemma says $V$ is a ISS Lyapunov function for the auxiliary 
system (\ref{LP3}). See \cite[Proof of Lemma 2.14]{SonWan1995} 
for instance.
Consider now the scaled error coordinates
$\varepsilon =\L^{-1}(\hat{z}-z)$.
Straightforward computations
from (\ref{eq_Syste})
give the error system
$$
\frac{1}{L}\,  \dot \varepsilon  \in S_m \varepsilon + \D_L  + \KR(\varepsilon_1+v)
$$
with $\D_L=\L^{-1}\delta$. 
Since $\Phi$ satisfies $\P(\alpha,\cc)$,
with (\ref{eq_defalphaHomo}) and $\frac{r_{i+1}}{r_j}\leq 1$,
we obtain, for all $L\geq 1$
\begin{align*}
|\D_{L,i}| 
&\leq \frac{\cc}{L} \sum_{j=1}^i L^{(j-1)\frac{r_{i+1}}{r_j}-i+1}{|\varepsilon_j|^{\frac{r_{i+1}}{r_j}}} + \frac{|\Delta w_i|}{L^i}  \ ,\\
&\leq \frac{\cc}{L} \sum_{j=1}^i{|\varepsilon_j|^{\frac{r_{i+1}}{r_j}}} + \frac{|\Delta w_i|}{L^i}  \ ,\\
&\leq \frac{c}{L} V(\varepsilon)^{\frac{r_{i+1}}{d_V}} +  \frac{|\Delta w_i|}{L^i} \ ,
\end{align*}
where $c$ is a positive real number obtained from Lemma \ref{lem2}
in Appendix \ref{app_lem}. With Lemma \ref{Lemm_Rob}, where $\bardelta _i$ plays the role of 
$\D_{L,i}$, $\barv$ the role of $v$ and $\bare$ the role of 
$\varepsilon $, we obtain that, by 
picking $L^*$ sufficiently large such that $\frac{c}{L^*}\leq \frac{c_\delta}{2}$, we have, for all $L>L^*$,
\\[0.7em]$\displaystyle 
\text{if }\left\{
\begin{array}{rl}
\displaystyle \frac{ |\Delta w_i| }{L^i}&\leq \displaystyle \frac{c_\delta }{4}  V(\varepsilon)^\frac{r_{i+1}}{d_V}\ ,\ \forall i\\
\displaystyle |v|&\leq \displaystyle c_v V(\varepsilon)^\frac{r_{1}}{d_V}
\end{array}\right.
$\refstepcounter{equation}\label{eq_ISSFinal}\hfill$(\theequation)$\\[0.5em]\null \hfill $\displaystyle 
\Rightarrow\frac{1}{L}\,  \max\left\{\frac{\partial V}{\partial e}(\varepsilon) \dot \varepsilon\right\} \leq -\lambda V(\varepsilon)^\frac{d_V+d_0}{d_V}\ .
$\\[0.7em]
Now, when evaluated along a solution, $\varepsilon $ gives rise to an absolutely continuous 
function $t\mapsto \varepsilon (t)$. Similarly 
the function defined by
$
t\mapsto \nu (t)=V(\varepsilon(t))
$
is absolutely continuous. It follows that its time derivative is defined for almost all $t$ and, according to \cite[p174]{Smi},
\eqref{eq_ISSFinal} implies, for almost all $t$,
\\[0.7em]$\displaystyle 
\text{if }\left\{
\begin{array}{rl}
\displaystyle \frac{ |\Delta w_i| }{L^i}&\leq \displaystyle \frac{c_\delta }{4}  \nu (t)^\frac{r_{i+1}}{d_V}\ ,\ \forall i\\
\displaystyle |v|&\leq \displaystyle c_v \nu (t)^\frac{r_{1}}{d_V}
\end{array}\right.
$\refstepcounter{equation}\label{LP7}\hfill$(\theequation)$\\\null \hfill $\displaystyle 
\Rightarrow
\frac{1}{L}\dot{\nu }(t) \leq -\lambda
\nu (t)^\frac{d_V+d_0}{d_V}\ .\quad \null 
$\\[0.7em]

Here two cases have to be distinguished.
\begin{enumerate}
\item 
If $d_0$ is in $]-1,0]$, with Lemma \ref{lem_LyapIss} in Appendix \ref{app_lem} (see also \cite{SonWan1995}), we
get the existence of a class $\mathcal{KL}$ function $\beta_V$ such that\footnote{according to Lemma  \ref{lem_LyapIss},
$\beta_V(s,t)=\max\{0,s^{\frac{-d_0}{d_V}}-t\}^{\frac{d_V}{-d_0}}$}
\\[0.7em]$\displaystyle 
V(\varepsilon(t)) \leq \max_{i\in[1,m]}
\left\{ 
\vrule height 1.2em depth 1.2em width 0pt
\right.\!
\beta_V(V(\varepsilon(0)),\lambda L t),
$\hfill \null \\\null \hfill $\displaystyle \sup_{s\in[0,t]}\left\{ 
\left(\frac{4 |\Delta w_i(s)| }{L^i c_\delta}\right)^\frac{d_V}{r_{i+1}} , \frac{|v(s)|^\frac{d_V}{r_{1}}}{c_v}\right\} 
\!
\left.
\vrule height 1.2em depth 1.2em width 0pt
\right\}\ .
$\\[0.7em]
The result holds since with Lemma \ref{lem2} there exist a positive real number $c_1$ such that
$$
\left|\frac{e_i}{L^{i-1}}\right|\leq c_1 V(\epsilon)^\frac{r_i}{d_V}\ .
$$
Moreover, when $v(t)=\Delta w_j(t)=\w_j=0$ for $j=1,\dots,m$,
(\ref{LP7}) implies finite time convergence
in the case in which $d_0<0$.\\
\item If $d_0=-1$, then $r_{m+1}=0$. We choose $L^*$ sufficiently 
large to satisfy
$$
\frac{ \bar w_m}{(L^*)^m}\leq \frac{c_\delta }{4}\ . 
$$
We obtain that the first condition in  (\ref{LP7}) is satisfied for 
$i=m$  when $L\geq L^*$.
With Lemma \ref{lem_LyapIss} in Appendix \ref{app_lem} (see also \cite{SonWan1995}),
the implication (\ref{LP7})
implies the existence of a class $\mathcal{KL}$ function $\beta_V$ such that\footnotemark[\thefootnote]
\\[0.7em]$\displaystyle 
V(\varepsilon(t)) \leq \max_{i\in[1,m-1]}
\left\{ 
\vrule height 1.2em depth 1.2em width 0pt
\right.\!
\beta_V(V(\varepsilon(0)),\lambda L t), 
$\hfill \null \\\null \hfill $\displaystyle \sup_{s\in[0,t]}\left\{ 
\left(\frac{4 |\Delta w_i(s)| }{L^i c_\delta}\right)^\frac{d_V}{r_{i+1}} , \frac{|v(s)|^\frac{d_V}{r_{1}}}{c_v}\right\} 
\!
\left.
\vrule height 1.2em depth 1.2em width 0pt
\right\}\ .
$\\[0.7em]
And the result holds as in the previous case.
\end{enumerate}
\end{pf}

\section{Cascade of homogeneous observers}
\label{sec_PR}
When we cannot find
$d_0$ in $[-1,0]$ and $\cc$ such that the nonlinearities satisfy $\P(\alpha,\cc)$, with $\alpha $ defined in (\ref{eq_defalphaHomo}),
we may lose
 the convergence of observer \eqref{eq_ObsHomo}, or the possibility of making the
final error arbitrarily small. In such a bad case, we can
still take advantage of the fact that, for $\alpha$ verifying
(\ref{eq_defalphaHomo}) with $d_0=-1$, $\P(\alpha,\cc)$ does not
impose any restriction besides boundedness of the last functions
$\Phi_m$ (see Table \ref{tab_holderPower}).

From the remark that observer \eqref{eq_ObsHomo}
\begin{enumerate}
	\item
	can be used for 
	the system
	$$\begin{array}{rcl}
	\dot{z}_1&=&z_2 + \psi_1(t)\\[0.321em]
	&\vdots&  \\[0.321em]
	\dot{z}_{k-1}&=&z_{k} + \psi_{k-1}(t)
	\\[0.321em]
	\dot{z}_{k}&=&\varphi_k(t) 
	\end{array}
	$$
	provided 
	the functions $\psi_i$ are known and the function $\varphi_k$ is unknown 
	but bounded, with known bound.
	\\[-0.5em]\null 
	\item
	gives estimates of the $z_i$'s in finite time,
\end{enumerate}
we see that it can be used as a preliminary step to deal with the 
system
\\[0.7em]
\vbox{\hsize=\linewidth\noindent
	\null \hfill $\begin{array}{rcl}
	\dot{z}_1&=&z_2 + \psi_1(t)\\[0.321em]
	&\vdots&  \\[0.321em]
	\dot{z}_{k-1}&=&z_{k} + \psi_{k-1}(t)
	\\[0.321em]
	\dot{z}_{k}&=&z_{k+1}+\Phi_{k}(u,z_1,\ldots,z_{k})  
	\\[0.321em]
	\dot{z}_{k+1}&=&\varphi_{k+1}(u,z_1,\ldots,z_{k+1}) 
	\end{array}$\hfill \null }\\[0.7em]
Indeed, thanks to the above observer we know in finite time the values 
of $z_1,\ldots,z_{k}$, so that the function $\Phi_k(u,z_1,\ldots,z_{k}) $ becomes a known signal $\psi_k(t)$.

From this, we can propose the following observer made of a cascade of 
homogeneous observers~:
\begin{equation}
\label{eq_ObsHomo_PR}
\renewcommand{\arraystretch}{1.4}
\begin{array}{@{}r@{\; }c@{\; }l@{}}
\dot{\hat{z}}_{11}&\in& \hat{w}_1-L_1 \, k_{11}\,  \sign(\hat z_{11}-y)
\\[-0.7em]
\multicolumn{3}{@{}c@{}}{\dotfill}\\
\dot{\hat{z}}_{21}&=&{\hat z}_{22} +  \Phi_1  (u,{\hat z}_{11})+\hat{w}_1 - L_2 \, k_{21}\,  \power{\hat z_{21}-y}{\frac{1}{2}}
\\
\dot{\hat{z}}_{22}&\in& \hat{w}_2 - L_2^2 \, k_{22}\,  \sign(\hat z_{21}-y)
\\[-0.7em]
\multicolumn{3}{@{}c@{}}{\dotfill}\\
&
\vdots
\\[-0.7em]
\multicolumn{3}{@{}c@{}}{\dotfill}\\
\dot{\hat{z}}_{m1}&=&{\hat z}_{m2}  +  \Phi_1  (u,{\hat z}_{11}) 
\\
\multicolumn{3}{r@{}}{+\hat{w}_1 - 
L_m \, k_{m1}\,  \power{\hat z_{m1}-y}{\frac{m-1}{m}}
}
\\[-0.4em]
& \vdots
\\[-0.4em]
\dot{\hat{z}}_{m(m-1)}&=&{\hat z}_{mm}  +
\Phi_{m-1}  (u,\hat z_{(m-1)1},\ldots,\hat z_{(m-1)(m-1)}) 
\quad \null 
\\
\multicolumn{3}{r@{}}{
  +\hat{w}_{m-1}	- L_m^{m-1} \, k_{m(m-1)}\,  \power{\hat z_{m1}-y}{\frac{1}{m}}
}
\\
\dot{\hat{z}}_{mm}&\in& \hat{w}_m
- L_m^m \, k_{mm}\,  \sign(\hat z_{m1}-y)
\end{array}
\end{equation}
where the $k_{ij}$ and $L_i$ are positive real numbers to be tuned.

As a direct consequence of Proposition \ref{prop_HM_pert} and following the same steps as in the proof of Proposition \ref{prop_HG_PR}, we have

\startrevision
\begin{proposition}
\label{prop_cascade}
 Assume $\Phi$ is continuous. 
 For any positive real numbers $\bz$, $\bu$ $\bw$,   
we can find positive real numbers  $k_{ij}$ and  $L_i$, two class $\mathcal K$ functions $\gamma_1$ and $\gamma_2$
and a class $\mathcal {KL}$ function $\beta$ such that, for all
locally  bounded time function $(u,v,w,\hat{w})$,
and all  $(z,\hat z)$ in
$\RR^m\times\RR^m$,
the observer \eqref{eq_ObsHomo_PR} admits absolutely continuous solutions
$\left(\vechatZsol_1(\hat z,z,t;u,v,w,\hat{w}),...,\vechatZsol_m(\hat z,z,t;u,v,w,\hat{w})\right)$ which are 
defined on $\RR_+$
and for any such solution we have for all $i$ in $\{1,...,m\}$ and for all $t$ such that $|Z(z,s;u,w)|\leq \bz$, $|u(s)|\leq \bu$ and $|\Delta w(s)|\leq \bw$ for all $0\leq s\leq t$:
\\[0.7em]$\displaystyle 
|\mathbf{\hat Z} _{{i}}(t)-\mathbf{Z}_i(t)| \leq \max\Big\{\beta(|z-\hat z|,t),
$\hfill \null \\\null \hfill $\displaystyle 
\sup_{\stackrel{1\leq j\leq i-1}{s\in[t_0,t]}} \left\{\gamma_1(|v(s)|), \gamma_2(|\Delta w_j(s)|)\right\}\Big\}. 
$\\[0.7em]
where $\vechatZsol_i$ is the state of the $i$th block (see Notation \ref{notation3}) and we have used the abbreviation
$\vechatZsol_i(t)=\vechatZsol_i(\hat z,z,t;u,v,w,\hat{w})$ and $\vecZsol_i(t)=\vecZsol_i(z,t;u,w)$.\\
Moreover, when  $v(t)=\Delta w_j=0$,	 there exists $T$ such that
$\mathbf{\hat Z} _{{i}}(\hat z,z,t)=\mathbf{Z}_i(z,t)$ for all $t\geq T$.
\end{proposition}
\stoprevision

This observer is an extension of the cascaded high gain observer \eqref{eqObsHG_PR} presented in Section
\ref{sec_HG_PR}.  The use of homogeneity enables here to obtain convergence without demanding anything but
the knowledge of a bound on the input and on the system solution.  A drawback of a
cascade of observers is that it gives an observer with dimension $\frac{m(m+1)}{2}$ in general.  However, as
seen in Section \ref{sec_HG_PR}, it may be possible to reduce this dimension since, for each new block, one
may increase the dimension by more than one, when the corresponding added functions $\Phi_i$ satisfy
$\P(\alpha,\cc)$ m for some $\alpha$ verifying (\ref{eq_defalphaHomo}) with $d_0=-1$
and for some $\cc$.

Finally, note that the result of Proposition \ref{prop_cascade} does 
not mean that the observer is ISS with respect to $\Delta w$. Indeed, 
$\Delta w$ must be bounded to obtain this ISS-like inequality : the 
system is ISS with restrictions.
Again, we believe that this problem could be solved employing homogeneous in the bi-limit observer as in \cite{AndPraAstSIAM}.

\startrevision
\section{Relaxing the assumptions marked with $\hypo$}
\label{sec_relaxHypo}
First, if System \eqref{eqCanForm} is not complete, every ISS inequalities still holds for any solution $Z(z,t;u,w)$ but only on $[0,T(z)[$ where $T(z)$ is its maximal time of existence. 

The global aspect of boundedness, H\"{o}der, $\P(\alpha,\cc)$, \ldots, can be relaxed as follows.
Let $U$ be bounded and  let $\M$ be a given compact set.
We define $\hat{\Phi}$, to be used  instead of $\Phi$ in the observers, as
\begin{equation}
\hat{\Phi}_i(u,z_1,...,z_i)=\sat(\Phi_i(u,z_1,...,z_i),\bar\Phi_i)
\label{eq_phiHat}
\end{equation}
where $\bar{\Phi}_i=\max_{u\in U,z\in\M}(\Phi_i(u,z_1,...,z_i))$ and the saturation function is defined on $\RR$ by
$$
\sat(x,M)=\max(\min(x,M),-M) \ . 
$$
It can be shown that, for any compact set $\tilde{\M}$ strictly contained in $\M$, there exists $\tilde{\cc}$
such that \eqref{eqAss_pertNew} holds for $\hat{\Phi}$ for all $(\vecZ_{{a}},\vecZ_{{b}})$ in
$\RR^m\times\tilde \M$.  Then, since $\hat{\Phi}=\Phi$ on $\tilde{\M}$,
we can modify the assumptions
\begin{list}{}{%
		\parskip 0pt plus 0pt minus 0pt%
		\topsep 0.5ex plus 0pt minus 0pt%
		\parsep 0pt plus 0pt minus 0pt%
		\partopsep 0pt plus 0pt minus 0pt%
		\itemsep 0pt plus 0pt minus 0pt%
		\settowidth{\labelwidth}{-}%
		\setlength{\labelsep}{0.5em}%
		\setlength{\leftmargin}{\labelwidth}%
		\addtolength{\leftmargin}{\labelsep}%
	}
	\item[-] in Proposition \ref{prop_HG_Nom}, so that \eqref{eq_LipscPerturbed} holds only on the 
	compact set $\M$;
	\item[-] in Propositions \ref{prop_HG} and \ref{prop_HM_pert}, so that $\Phi$ verifies 
	$\P(\alpha,\cc)$ only on the compact set $\M$;
	\item[-] Propositions \ref{prop_HG_PR} and \ref{prop_cascade} remain unchanged. 
\end{list}
\vskip -1.3em
In this case, the results hold for the particular system solutions $Z(z,t;u,w)$ which are in the compact set $\tilde{\M}$ for $t$ in
$[0,T(z)[$. Precisely, for these solutions, the bounds  on $\hat Z_i(t)-Z_i(t)$ given in these Propositions hold for all $t$ in $[0,T(z)[$.

Note also that if $\P(\alpha,\cc)$ holds on a compact set, then for any $\tilde{\alpha}$ such that $\tilde{\alpha}_{ij}\leq \alpha_{ij}$ for all $(i,j)$, there exists $\tilde{a}$ such that $\P(\tilde\alpha,\tilde\cc)$ also holds on this compact set. It follows that the constraints given by \eqref{eq_defalphaHomo} or Table \ref{tab_holderPower} in Proposition \ref{prop_HM_pert} can be relaxed to 
$
\alpha_{ij}  \geq 
\frac{1-d_0(m-i-1)}{1-d_0(m-j)}$, 
and the less restrictive conditions one may ask for are obtained for $d_0=-1$. 

Finally, in Propositions \ref{prop_HG_Nom}, \ref{prop_HG} and \ref{prop_HM_pert}, it is possible to consider the case where $\Phi$ depends also on time as long as any assumption made on $\Phi$ is satisfied uniformly with respect to time.
\stoprevision

\section{Example}
\label{sec_ex}

Consider the system
\begin{equation}
\label{eq_ex2}
\dot{x}_1 \:=\:  x_2 \; ,\  
\dot{x}_2 \:=\:  -x_1 + x_3^5 x_1 \; ,\  
\dot{x}_3 \:=\:  -x_1x_2 +u \; ,\  
y \:=\:  x_1 
\end{equation}
with $u$ as input.
It would lead us too far from the
main subject of this article to study here the solutions behavior
of this system. We note however that, when $u$ is zero, they evolve in 
the $2$-dimensional surface
$\{x\in \RR^3: 3x_1^2+3x_2^2+x_3^6=c^6\}$ which is 
diffeomorphic\footnote{%
	A diffeomorphism from the unit sphere to the set is $x\mapsto x \rho 
	(x)$ where $\rho $ is the unique positive solution (hint: $x_3\leq 1$) of
	$
	\rho ^6 x_3^6\;+\; 3\rho ^2 (1-x_3^2)\;-\; 1\;=\; 0
	$
}
to 
the sphere $\SS^2$.
Thanks to 
Poincar\'{e}-Bendixon theory, we know the solutions are 
periodic and circling the  unstable
 equilibria $(x_1=x_2=0,x_3=\pm c)$.
So we hope for the 
existence of solutions remaining in the compact set
$$
\C_{r,\epsilon}\;=\; \left\{x\in\RR^3:\,  x_1^2+x_2^2 \geq \epsilon \ , \
3x_1^2+3x_2^2+x_3^6 \leq r\right\}
$$
for instance when $u$ is a small 
periodic time function, except maybe for pairs of 
input $u$ and
initial condition $(x_1,x_2,x_3)$  for which resonance could occur. 
Moreover, due to their periodic behavior, such solutions are likely to have their $x_3$ component recurrently crossing zero.

\subsection{Uniform and differential observability}

On $\Ouvs= \left\{x\in\RR^3:\,  x_1^2+x_2^2 \neq 0\right\}$, and whatever $u$ is, the knowledge of the 
function $t\mapsto y(t) = X_1(x,t)$ and therefore of its three first derivatives
\begin{eqnarray*}
	\dot{y} &=& x_2 \\
	\ddot{y} &=& -x_1+	x_3^5	x_1 \\
	\dddot{y}&=&-x_2-	5x_3^4	x_1^2x_2 + x_3^5x_2 +	5x_3^4	x_1u
\end{eqnarray*}   
gives us $x_1$, $x_2$ and $x_3$. Thus, System  (\ref{eq_ex2}) is uniformly observable on $\Ouvs$. 
Besides, the function
$$
\mathbf{H}_4(x) = \left(\begin{array}{c}
x_1 \\ x_2 \\ -x_1+x_3^5x_1 \\ -x_2-5x_3^4x_1^2x_2 + x_3^5x_2
\end{array}\right)
$$
is injective on $\Ouvs$ and admits the left following inverse, defined on $\left\{z\in\RR^4:\,  z_1^2+z_2^2 \neq 0\right\} $, is:
$$
\mathbf{H}^{-1}_4(z) = \left(\begin{array}{c}
z_1 \\z_2 \\
\left(\frac{(z_3+z_1)z_1 + \left[(z_4+z_2)+ 
	3|(z_3+z_1)\power{z_1}{\frac{3}{2}}|^{\frac{4}{5}}z_2\right]z_2}{z_1^2+z_2^2}\right)^{\frac{1}{5}}
\end{array}\right)
$$
However, $\mathbf{H}_4$ is not an immersion because of a singularity of its Jacobian 
at $x_3=0$. So the system is 
differentially observable of order $4$ on $\Ouvs$ but not strongly. 
According to  \cite{BerPraAndAUT20161}, it admits a triangular canonical form of dimension $4$ but with 
functions $\Phi$ maybe non Lipschitz. 

\subsection{Triangular form and property $\P(d_0,c,0)$}

The triangular canonical form of dimension $4$ mentioned above is
\begin{equation}
	\label{LP13}\begin{array}{rcl}
		\dot{z}_1 &=& z_2 \\
		\dot{z}_2 &=& z_3 \\
		\dot{z}_3 &=& 
		z_4+\Phi_3(u,z_1,z_2,z_3) 
		\\
		\dot{z}_4 &=& 
		\Phi_4(u,z)\\
		y &=& z_1 .
	\end{array}
\end{equation}
where $\Phi_3(u,z_1,z_2,z_3)=5u |z_3+z_1|^{\frac{4}{5}}
\power{z_1}{\frac{1}{5}}$ and $\Phi_4$ is a continuous non-Lipschitz
function the expressions of which is complex, fortunately with no
interest here.  The function $\Phi_3$ is not Lipschitz at the points
on the hyperplanes $z_3=-z_1$ and $z_1=0$ (image by $\mathbf{H}_3$ of
points where $x_3=0$ or $x_1=0$) known to be visited possibly
recurrently along solutions.  This example thus falls precisely into
the scope of the paper.

The function $\Phi_4$ is continuous and therefore bounded on any compact set including $\textbf{H}_4(\bar \C_{r,\epsilon})$. 
Besides, for
$\vechatZ_{{3}}$ and $\vecZ_{{3}}$ in a compact set including 
$\textbf{H}_3(\bar \C_{r,\epsilon})$, there exist $c_1$ and $c_3$ such that
\\[0.7em]$\displaystyle 
|\Phi_3(u,\hat{z}_1,\hat{z}_2,\hat{z}_3)-\Phi_3(u,z_1,z_2,z_3)| 
$\hfill \null \\\null \hfill $\displaystyle  \leq c_1 u |\hat{z}_1-z_1|^{\frac{1}{5}}
+c_3 u |\hat{z}_3-z_3|^{\frac{4}{5}}
\  .
$\\[0.7em]
This implies that $\Phi_3$ is H\"older with order $\frac{1}{5}$. 

Hence the nonlinearities $\Phi_3$ and $\Phi_4$ verify the conditions 
of Table \ref{tab_holderPowerPract}. This implies that for $L$ 
sufficiently large, convergence with an arbitrary small error  can be achieved with the high gain observer \eqref{eq_ObsHG} . However, $\Phi_3$ does not verify the conditions of Table \ref{tab_holderPower}.
Thus, there is no theoretical guarantee that  the homogeneous observer \eqref{eq_ObsHomo} with $d_0=-1$ will  provide exact convergence.

\subsection{An observer of dimension 4 ?}

\begin{center}
	\begin{figure}
		\includegraphics[width=\columnwidth]{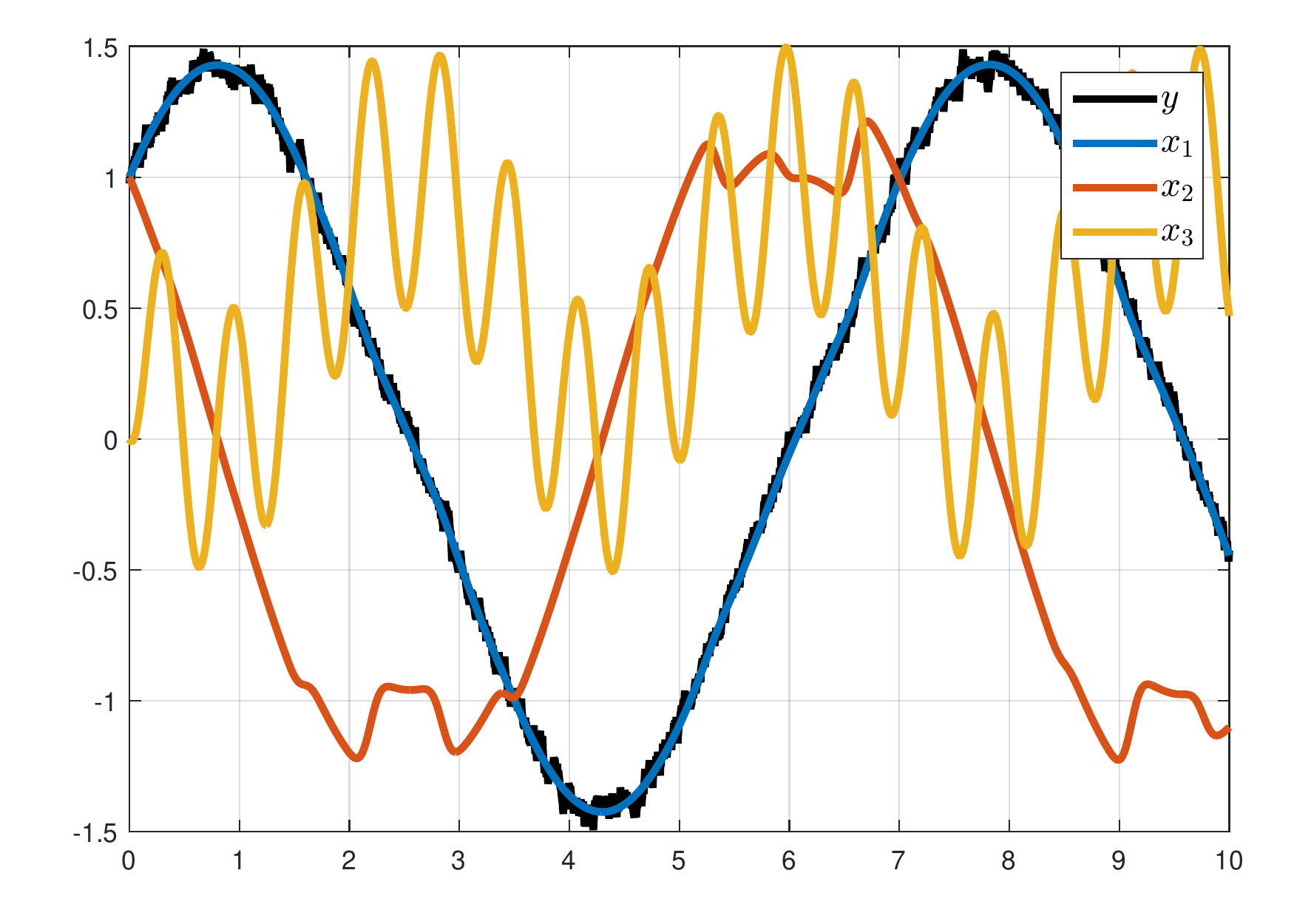}
		\caption{Trajectory of System \eqref{eq_ex2}, with the noised measurement $y$.}
		\label{fig_exTraj}
	\end{figure}
\end{center}

We consider the solution to system \eqref{eq_ex2} with initial condition $x=(1,1,0)$ and $u=5\sin(10 \, t)$. This solution is periodic and regularly crosses the Lipschitzness singularities 
$x_3=0$ or $x_1=0$,
as illustrated in Figure \ref{fig_exTraj}. In the following, we use the same noised measurement $y$, shown on Figure \ref{fig_exTraj}, in every simulation with noise. It is a filtered gaussian noise with standard deviation $\sigma=0.03$ and $1$st order filtering parameter $a=50$.

We first implement a high gain observer of dimension 4, in the absence of noise,  initialized 
at $\hat{x}=(0.1,0.1,0)$, and with the gains $k_1=14$, $k_2=99$, $k_3=408$, $k_4=833$.
As an illustration of Proposition \ref{prop_HG}, the convergence with
an arbitrary small error is achieved and is illustrated in Table
\ref{tab_exHG}.  However, we observe that the decrease of the errors,
especially for $e_{z,4}$, is very slow compared to the increase of the
peaking and a very high gain is needed to obtain "acceptable" final
errors.  In presence of noise, the tradeoff between final error and
noise amplification becomes very difficult : with the noised
measurement of Figure \ref{fig_exTraj}, the smallest final error
$e_{z,4}$ is $200$, achieved for $L=2$.  Of course, there might
exist a choice of the gains $k_i$ giving better results.  But overall
a high gain observer may not be a systematic solution in practice for
non-Lipschitz triangular systems, especially when the solution
regularly crosses the Lipschitz-singularities.
\begin{table}[h]
	\begin{center}
		\begin{tabular}{|c|c|c|c|c|c|}
			\hline
		L	& $e_{z,1}$ & $e_{z,2}$ & $e_{z,3}$ & $e_{z,4}$ & Peaking \\
		\hline
		2  &   0.15 &  4  & 60  & 200 &  300 \\
		\hline
		5 & $6 \,.\, 10^{-4}$ & 0.04 & 1.5 & 30 & 4000  \\
		\hline
		8 & $5 \,.\, 10^{-5}$ & $4 \,.\, 10^{-3}$ & 0.25& 7 & $1.5 \,.\, 10^{4}$ \\
		\hline
		10 & $8 \,.\, 10^{-6}$ & $1\,.\, 10^{-3}$ & 0.1 &4 & $3.5 \,.\, 10^{4}$ \\
		\hline
		15 & $1.5 \,.\, 10^{-6}$ & $3\,.\, 10^{-4}$ & 0.03 & 2 & $1.2 \,.\, 10^{5}$\\
		\hline
		\end{tabular}
	\end{center}
\caption{Decrease of the final error in the $z$-coordinates ($e_{z,i}=\hat{z}_i-z_i$) with the gain $L$, with a high gain observer and in the absence of noise.}
\label{tab_exHG}
\end{table}

Let us now implement an homogeneous observer of dimension $4$ with an explicit Euler method with fixed measurement and integration steps equaling $10^{-5}$, and with the Matlab sign function. The degree is $d_0=-1$, and the gains are chosen according to \cite{Lev2005}, i-e $k_1=5$, $k_2=8.77$, $k_3=4.44$, $k_4=1.1$. For a gain $L=3$, the convergence is achieved with a final error of $8 \,.\,10^{-4}$ on $z_4$, even though the H\"older restriction of Proposition \ref{prop_HM_pert} is a priori not satisfied around $z_1=0$. Unfortunately, the final errors are heavily impacted in presence of noise, as illustrated in Table \ref{tab_exHM}.
This may also come from a lack of ISS property.  
Notice that the amplification of the noise by the gain $L$ is not as rapid as expected from the bound in Proposition \ref{prop_HM_pert}. The final errors remain nonetheless too large, although, once again, we did not optimize our choice of gains $k_i$.
\begin{table}[h]
	\begin{center}
		\begin{tabular}{|c|c|c|c|c|c|}
			\hline
			L	& $e_{z,1}$ & $e_{z,2}$ & $e_{z,3}$ & $e_{z,4}$ \\
			\hline
			2.5  &   0.15 & 3.5  & 30  & 18 \\
			\hline
			3  &   0.15 & 3  & 35  & 25 \\
			\hline
			4 & 0.1 & 2 & 25 & 50  \\
			\hline
			5 & 0.1 & 2 & 30& 80 \\
			\hline
			6 & 0.1 & 2 & 35 &120 \\
			\hline
		\end{tabular}
	\end{center}
	\caption{Final errors in the $z$-coordinates given by a homogeneous observer of degree $-1$ in presence of  noise.}
	\label{tab_exHM}
\end{table} 

\subsection{Cascaded observers}

In the absence of noise, the cascaded observers presented in Sections \ref{sec_HG_PR} and \ref{sec_PR} give 
similar results to the corresponding observers in dimension 4, i-e  arbitrary small asymptotic error and finite time convergence respectively. However, they seem to provide better accuracies in presence of noise.  

In the case of a high gain cascade observer, the errors, although 
smaller than in the high gain observer of dimension 4, remain too large to consider it a viable solution. On the other hand, the homogeneous cascade observer :
$$
\renewcommand{\arraystretch}{1.4}
\begin{array}{@{}r@{\; }c@{\; }l@{}}
\dot{\hat{z}}_{11}&=&{\hat z}_{12} - L_1 \, k_{11}\,  \power{\hat z_{11}-y}{\frac{2}{3}}
\\
\dot{\hat{z}}_{12}&=&{\hat z}_{13} - L_1^2 \, k_{12}\,  \power{\hat z_{11}-y}{\frac{1}{3}}
\\
\dot{\hat{z}}_{13}&\in& - L_1^3 \, k_{13}\,  \sign(\hat z_{11}-y)
\\[-0.7em]
\multicolumn{3}{@{}c@{}}{\dotfill}\\
\dot{\hat{z}}_{21}&=&{\hat z}_{22} - L_2 \, k_{21}\,  \power{\hat z_{21}-y}{\frac{3}{4}}
\\
\dot{\hat{z}}_{22}&=&{\hat z}_{23} - L_2^2 \, k_{22}\,  \power{\hat z_{21}-y}{\frac{1}{2}}
\\
\dot{\hat{z}}_{23}&=& {\hat{z}}_{24} + \sat(\gz_3({\hat z}_{11},{\hat z}_{12},{\hat z}_{13})) u - L_2^3 \, k_{23}\, \power{\hat z_{21}-y}{\frac{1}{4}}
\\
\dot{\hat{z}}_{24}&\in& - L_2^4 \, k_{24}\,  \sign(\hat z_{21}-y)
\end{array}
$$
with the coefficients $k_{1j}$ chosen, according to \cite{Lev2005}, as $k_{11}=3$, $k_{12}=2.6$, $k_{13}=1.1$, and $k_{2j}$ as above, and with the gains $L_1=2.5$ and $L_2=3$, gives the following final errors :
$$
e_{z,11}= 0.05,\quad e_{z,12}= 0.4,\quad e_{z,13}= 2.5,\quad e_{z,24}= 12 
$$
Comparing to Table \ref{tab_exHM}, we see that implementing an intermediate homogeneous observer of dimension 3 enables to obtain much better estimates of the first three states $z_i$, which are then used in the nonlinearity of the second block, thus giving a better estimate of $z_4$.

Unfortunately, the presented results are still unsatisfactory in presence of noise, which leaves the question of the construction of robust observers for such systems unanswered.

\section{Conclusion}
\startrevision
To summarize the most important ideas, we provide in Table \ref{tab_compObs} a synthetic comparison of the four observers studied in this paper, in the usual case where the system state and the input are bounded.
\stoprevision

We have shown the convergence with an arbitrary small error of the classical high gain observer in presence of nonlinearities verifying some H\"older-like condition. The same result could probably be obtained for the high gain-like observer presented in \cite{Bes}. Also, for the case when this H\"older condition is not verified, we proposed a novel cascaded high gain observer. 
Under slightly more restrictive assumptions, we proved the convergence of an homogeneous observer and of its cascaded version with the help of an explicit Lyapunov function. 

\startrevision
Our numerical experience indicates however that to improve the performances in presence of
measurement noise, it is very difficult to tune the gains of both high gain and homogeneous observers,
although it is slightly simpler for the latter since smaller gains are sufficient to ensure convergence.
Simulations on our example suggest that the situation may be more favorable with the cascaded homogeneous
observer.
Our ISS bounds in this paper being far too conservative, it is necessary to carry out a finer study
if we want to optimally tune the gains of the observers. It may also be
appropriate to use on-line gain adaptation techniques since large gains should be necessary only
around the points where the nonlinearities are not Lipschitz. About these two aspects,
we refer the reader to the survey in 
\cite[Sections 3.2.2 and 3.2.3]{KhaPra} and the references therein.
\stoprevision

\begin{table*}[ht]
\startrevision
	\centering
	\begin{tabular}{|p{0.10\linewidth}|p{0.2\linewidth}|p{0.2\linewidth}|p{0.2\linewidth}|p{0.2\linewidth}|}
		\cline{2-5} \multicolumn{1}{c|}{} & High gain \eqref{eq_ObsHG} & High gain cascade \eqref{eqObsHG_PR} & Homogeneous \eqref{eq_ObsHomo} & Homogeneous cascade \eqref{eq_ObsHomo_PR} \\
		\hline
		 Assumption on $\gz_i$& H\"older with order greater than in Table \ref{tab_holderPowerPract} &Continuous  &H\"older with order greater than in \eqref{eq_defalphaHomo} or Table \ref{tab_holderPower} for $d_0=-1$  &Continuous \\[0.5em]
		 Convergence & Arbitrary small error & Arbitrary small error&Asymptotic convergence & Asymptotic convergence \\[0.5em]
		 Advantages& Easy choice of gains & No constraint on $\gz_i$  & Not necessarily large gains because convergence & No constraint on $\gz_i$,  convergence, apparently better in terms of noise \\
		 Drawbacks& Large gains necessary to obtain small error $\Rightarrow$ numerical problems (peaking) and sensitivity to noise &Same as for high gain, but also gains difficult to choose and large dimension & Implementation of the sign function if $d_0=-1$ (chatter etc) &Large dimension and a lot of gains to choose\\
		\hline
	\end{tabular}
		\label{tab_compObs}
		\caption{Comparison between observers when the system state and the input are bounded.}
\stoprevision
\end{table*}

\appendix

\section{Barbot et al's observer}
\label{app_Barbot}

The set valued map proposed in \cite{BarBouDje} to obtain an 
observer for a triangular canonical form where the functions are only locally 
bounded is defined as follows.
Given  $(\hat{z},y,u)$, 
$(v_1,\ldots,v_m)$ is in 
$F(\hat{z},y,u)$ 
if there exists $(\tilde{z}_2,...,\tilde{z}_m)$ in $\RR^{m-1}$ such that:
\begin{eqnarray*}
	v_1&=&\tilde{z}_2 + \gz_1(y) \, u  
	\\
	\tilde{z}_2&\in& \sat(\hat{z}_2) -k_1 \, S(y-\hat{z}_1) 
	\\
	&\vdots& 
	\\
	v_i&=&\tilde{z}_{i+1} + \gz_i(y,\tilde{z}_2,\ldots,\tilde{z}_i) \, u 
	\\
	\tilde{z}_{i+1}&\in&\sat(\hat{z}_{i+1}) -k_i \, S(\hat{z}_i-\tilde{z}_i) 
	\\
	&\vdots& 
	\\
	v_m&\in& \varphi_m(y,\tilde{z}_2,\ldots,\tilde{z}_m) 
	\\
	&& \quad+ \gz_m(y,\tilde{z}_2,\ldots,\tilde{z}_m) \, u -k_m \, S(\hat{z}_m-\tilde{z}_m)  
\end{eqnarray*}
where $\sat$ is some saturation function.

\section{Proof of Lemma \ref{Lemm_Rob}}\label{Sec_ProofLem1}
The proof is based on the following Proposition the proof of which is
given in the following section for the case $d_0=-1$ and can be found
for $d_0$ in $]-1,0]$ in \cite{AndPraAstSIAM}.  This proposition
establishes that for a chain of integrator it is possible to construct
homogeneous correction terms which provide an observer and that it is
possible to construct a smooth strict homogeneous Lyapunov function.
\begin{proposition}\label{Lemm_StricteLyap}
For all $d_0$ in $[-1,0]$, the function $V$ defined in (\ref{eq_defV}) is positive definite and there exists positive real numbers $k_1, \dots k_m$, $\ell_1, \dots\ell_m$, $\tilde \lambda$ such that for all $e$ in $\RR^m$,
the following  holds~:
\begin{equation}
\label{LP4}
\max\left\{\frac{\partial V}{\partial \bare}(\bare)
\left(S_m(\bare) +  \KR(\bare_1)\right)\right\} \; \leq\;  -\tilde \lambda V(\bare)^\frac{d_V+d_0}{d_V}.
\end{equation}
\end{proposition}

Let  $\tilde{\KR}(\bare_1,s)$ be the function defined as
$$
\left(\tilde{\KR}(\bare_1,s)\right)_i = \left(\KR(\bare_1)\right)_i\ ,\ i\in [1,m-1]\ ,
$$
and,
$$
\left(\tilde{\KR}(\bare_1,s)\right)_m = \left\{\begin{array}{l}k_ms\ ,\ \text{ when } d_0=-1\\
\left(\KR(\bare_1)\right)_m\ ,\ \text{ when } d_0>-1
\end{array}\right. \ .
$$
Note that $\tilde{\KR}$ is a continuous (single) real-valued function which satisfies for all $\bare_1$ in $\RR$
$$
\KR(\bare_1)=\{\tilde{\KR}(\bare_1,s) \, , \quad s\in \sign(\bare_1)\}\ .
$$
Consider also the functions
\begin{multline*}
\tilde \eta(\bare,\bardelta,\barv,s) = \frac{\partial V}{\partial 
\bare}(\bare) (S_m(\bare) +  \bardelta + \tilde \KR(\bare_1+\barv,s)) \\+ \frac{\tilde \lambda}{2} V(\bare)^\frac{d_V+d_0}{d_V},
\end{multline*}
and
$$
\gamma(\bardelta,v) = \sum_{i=1}^m |\bardelta_i|^\frac{d_V +d_0}{r_{i+1}} + |\barv|^\frac{d_V+d_0}{r_1}.
$$
With (\ref{LP4}), we invoke Lemma \ref{lem2} to get the existence of 
a positive real number $c_1$ satisfying~:
\\[0.7em]$\displaystyle 
\frac{\partial V}{\partial \bare}(\bare) (S_m(\bare) +  \Delta + \tilde \KR(\bare_1+\barv,s)) 
$\hfill \null \\\null \hfill $\displaystyle \leq - \frac{\tilde \lambda}{2} V(\bare)^\frac{d_V+d_0}{d_V} 
+ c_1\sum_{i=1}^m \bardelta_i^\frac{d_V +d_0}{r_{i+1}} + c_1 |\barv|^\frac{d_V+d_0}{r_1}\ .
$\\[0.7em]
This can be rewritten,
\begin{multline*}
\frac{\partial V}{\partial \bare}(\bare) (S_m(\bare) +  \bardelta + \tilde \KR(\bare_1+v,s)) 
\leq -\frac{\tilde \lambda}{2(m+2)}  V(\bare)^\frac{d_V+d_0}{d_V}\\ 
+ \sum_{i=1}^m \left (c_1|\bardelta_i|^\frac{d_V +d_0}{r_{i+1}}-\frac{\tilde \lambda}{2(m+2)} V(\bare)^\frac{d_V+d_0}{d_V}\right) \\
+ c_1 |\barv|^\frac{d_V+d_0}{r_1}-\frac{\tilde \lambda}{2(m+2)} V(\bare)^\frac{d_V+d_0}{d_V}.
\end{multline*}
Consequently, the result holds with $\lambda = \frac{\tilde \lambda}{2(m+2)}$, $c_\delta = c_v = \left(\frac{\alpha}{c_1}\right)^\frac{r_1}{d_V+d_0}$.

\section{Proof of Proposition \ref{Lemm_StricteLyap} when $d_0=-1$}

In this section, we denote $\efin_i=(e_i,...,e_m)$.
Let $d_V$ be an integer such that
$d_V> 2m-1$  and the functions $\KR_i$ recursively defined by :
$$
\KR_m(e_m)=-\power{ e_m}{0} \quad ,\quad 
\KR_i(e_i)=\left(\begin{array}{c}
-\power{\ell_i e_i}{\frac{r_{i+1}}{r_i}} \\
\KR_{i+1}\left(\power{\ell_i e_i}{\frac{r_{i+1}}{r_i}}\right)
\end{array}
\right)\ .
$$
Let $V_m(e_m)=\frac{|e_m|^{d_V}}{d_V}$ and for all $i$ in $\{1,\dots,m-1\}$, let also $\bar V_i:\RR^2\rightarrow \RR$ and $V_i:\RR^{n-i+1}\rightarrow\RR$ be the functions defined by
\begin{align*}
\bar V_i(\nu,e_{i+1})&=\int_{\power{e_{i+1}}{\frac{r_i}{r_{i+1}}}}^{\nu} {\power{x}{\frac{d_V-r_i}{r_i}}- \power{e_{i+1}}{\frac{d_V-r_i}{r_{i+1}}} dx}\ ,\\
V_i(\efin_i) &= \sum_{j=m-1}^i \bar V_j(\ell_j e_j,e_{j+1}) + V_m(e_m)\ .
\end{align*}
With these definitions, the Lyapunov function $V$ defined in (\ref{eq_defV}) is simply $V(e)=V_1(e)$ and
the homogeneous vector field 
 $\KR(e_1)=\KR_1(e_1)$ with
$$
k_i = \ell_i^{\frac{r{i+1}}{r_i}} \ell_{i-1}^{\frac{r_{i+1}}{r_{i-1}}}\, ... \, \ell_2^{\frac{r_{i+1}}{r_2}} \ell_1^{\frac{r_{i+1}}{r_1}} \ .
$$
Note that the $j$th component of $\KR_i$ is homogeneous of degree 
$r_{j+1}=m-j$
and, for any $e_i$ in $\RR$, the set $\KR_{i}(e_{i})$ can be expressed as
$$
\KR_{i}(e_{i})=\{\tilde{\KR}_{i}(e_{i},s) \, , \quad s\in 
\sign(e_{i})\}\ ,
$$
where $\tilde{\KR}_{i}:\RR\times [-1,1]\to \RR$ is a continuous (single 
valued) function.

The proof of Proposition \ref{Lemm_StricteLyap} is made iteratively from $i=m$ toward $1$.
At each step, we show that $V_i$ is positive definite and we look for a positive real number $\ell_i$, such that 
  for all $\efin _i$ in $\RR^{n-i+1}$
\begin{multline}\label{eq_rec}
\max_{s\in\sign(e_i)} \left\{\frac{\partial V_{i}}{\partial \efin _{i}}(\efin _{i})(S_{m-i+1}\efin _{i} + \tilde{\KR}_{i}(e_{i},s))\right\} \\\leq - c_i V_i(\efin _i)^\frac{d_V-1}{d_V}\ ,
\end{multline}
where $c_i$ is a positive real number. The Proposition will be proved once we have shown that the former inequality holds for $i=1$.

\underline{Step $i=m$~:}
At this step, $\efin _m=e_m$.
Note that we have
\begin{align*}
\max_{s\in\sign(e_m)}\left\{ \frac{\partial V_{m}}{\partial \efin _{m}}(\efin _{m}) \tilde{\KR}_{m}(e_{m},s)\right\} &= - |\efin _m|^{d_V-1} \ ,\\
&=-c_m V_m(\efin _m)^\frac{d_V-1}{d_V}\ ,
\end{align*}
with $c_m = d_V^\frac{d_V-1}{d_V}$. Hence, equation (\ref{eq_rec}) holds for $i=m$.

\underline{Step $i=j$~:}
Assume $V_{j+1}$ is positive definite and assume there exists $(\ell_{j+1}, \dots, \ell_m)$ such that (\ref{eq_rec}) holds for $j=i-1$.
Note that the function $x\mapsto \power{x}{\frac{d_V-r_j}{r_j}}- \power{e_{i+1}}{\frac{d_V-r_j}{r_{j+1}}}$ is strictly increasing, 
is zero iff $x=\power{e_{j+1}}{\frac{r_j}{r_{j+1}}}$, and therefore has the same sign as $x-\power{e_{j+1}}{\frac{r_j}{r_{j+1}}}$. Thus, for any $e_{j+1}$ fixed in $\RR$, the function $\nu\mapsto \overline V_j (\nu,e_{j+1})$ is non negative and is zero only for $v=\power{e_{j+1}}{\frac{r_j}{r_{j+1}}}$. 
Thus, $\bar V_j$ is positive and we have
\begin{align*}
V_j(\efin _{j})=0 &\Leftrightarrow 
\left\{
\begin{array}{c}
V_{j+1}(\efin _{j+1})=0 \\
\overline{V}_j (\ell_j e_j,e_{j+1})=0
\end{array}
\right.\\
&\Leftrightarrow
\left\{
\begin{array}{l}
\efin _{j+1}=0 \\
\ell_j e_j=\power{e_{j+1}}{\frac{r_j}{r_{j+1}}}=0
\end{array}
\right.
\end{align*}
so that $V_j$ is positive definite. 

On another hand, let $\tilde V_j(\nu,\efin _{j+1})= V_{j+1}(\efin _{j+1}) + \bar V_j(\nu,e_{j+1})$
and let $T_1$ be the function defined
$$
T_1 (\nu,\efin _{j+1}) = \max_{s\in \sign(\nu)} \left\{\tilde T_1 (\nu,\efin _{j+1},s) \right\}
$$
with $\tilde T_1$ continuous and  defined by
\\[0.5em]\null \hfill $\displaystyle 
\tilde T_1 (\nu,\efin _{j+1},s) = 
\frac{\partial \tilde V_{j}}{\partial \efin _{j+1}}(\efin _{j+1})(S_{m-i-1}\efin _{i+1} 
+ \tilde{\KR}_{j+1}(\power{\nu}{\frac{r_{j+1}}{r_j}},s))
+\frac{c_{j+1}}{2}\tilde V_j(\nu,\efin _{j+1})^\frac{d_V-1}{d_V}\ .
$\hfill \null \\[0.5em]
Let also $T_2$ be the continuous real-valued function defined by
$$
T_2(v,\efin _{j+1}) = -\frac{\partial \tilde V_{j}}{\partial \nu}(\nu,\efin _{i+1})(   e_{j+1} - \power{\nu}{\frac{r_{j+1}}{r_j}})\ .
$$
Note that $T_1$ and $T_2$ are homogeneous with weight $r_j$ for $\nu$ and $r_i$ for $e_i$ and degree $d_V - 1$.
Besides, they verify the following two properties :
\begin{itemize}
	\item[-] for all $\efin _{j+1}$ in $\RR^{m-j}$, $\nu$ in $\RR$
\\[0.5em]\null \hfill $\displaystyle 
	T_2(\nu,\efin _{j+1})\geq 0
$\hfill \null \\[0.5em]
	(since $( \power{\nu}{\frac{r_{j+1}}{r_j}} -e_{j+1})$ and $(\power{\nu}{\frac{d_V-r_j}{r_j}}-  \power{e_{j+1}}{\frac{d_V-r_j}{r_{j+1}}})$ have the same sign)
	\item[-] for all $(\nu,\efin _{j+1})$ in $\RR^{m-j+1}\setminus\{0\}$, and $s$ in $\sign(\nu)$, we have the implication
\\[0.5em]\null \hfill $\displaystyle 
	T_2(\nu,\efin _{j+1})=0 \quad \Rightarrow \quad \tilde{T}_1(\nu, \efin _{j+1},s)<0
$\hfill \null \\[0.5em]
	since $T_2$ is zero only when $\power{\nu}{\frac{r_{j+1}}{r_j}}=e_{j+1}$ and
	\\[0.7em]$\displaystyle 
	\tilde{T}_1(\power{e_{j+1}}{\frac{r_{j+1}}{r_j}}, \efin _{j+1},s) =
	$\hfill \null \\\null \hfill $\displaystyle 
	\frac{\partial V_{j+1}}{\partial \efin _{j+1}}(\efin _{j+1})(S_{n-i}\efin _{j+1} + \tilde{\KR}_{j+1}(e_{j+1},s)) 
	$\hfill \null \\\null \hfill $\displaystyle
	+\frac{c_{j+1}}{2}V_{j+1}(\efin _{j+1})^\frac{d_V-1}{d_V}\leq -\frac{c_{j+1}}{2}V_{j+1}(\efin _{j+1})^\frac{d_V-1}{d_V}\ ,
	$\\[0.7em] 
	where we have employed (\ref{eq_rec}) for $i=j-1$.
\end{itemize}
Using Lemmas \ref{lem} in Appendix \ref{app_lem}, there exists $\ell_j$ such that
\\[0.5em]\null \hfill $\displaystyle 
T_1(\nu,\efin _{j+1}) - \ell_j T_2(\nu,\efin _{j+1}) \leq 0\ ,\ \forall \ (\nu,\efin _{j+1})\ .
$\hfill \null \\[0.5em]
Finally, note that
\\[0.5em]$\displaystyle 
\max_{s\in\sign(e_i)} \left\{\frac{\partial V_{j}}{\partial \efin _{j}}(\efin _{j})(S_{m-j+1}\efin _{j} + \tilde{\KR}_{j}(e_{j},s))\right\}  = 
$\hfill \null \\\null \hfill $\displaystyle 
 T_1(\ell_j e_j) - \ell_j T_2(\ell_je_j,\efin _{j+1}) - \frac{c_{j+1}}{2}V_{j}(\efin _{j})^\frac{d_V-1}{d_V}
$\\[0.5em]
Hence, (\ref{eq_rec}) holds for $i=j$.

\section{Technical lemmas}
\label{app_lem}

\begin{lem}
\label{lemBis}
Let $\eta$ be a continuous functions defined on $\RR^{n+1}$  and $f$ a continuous function defined on $\RR^n$. Let $\C$ be a compact subset of $\RR^n$. Assume that, for all $x$ in $\C$ and $s$ in $\sign(f(x))$,
\\[0.5em]\null \hfill $\displaystyle 
	\eta(x,s) < 0 \ .
$\hfill \null \\[0.5em]
Then, there exists $\alpha>0$ such that for all $x$ in $\C$ and $s$ in $\sign(f(x))$
\\[0.5em]\null \hfill $\displaystyle 
	\eta(x,s)  <-\alpha \ .
$\hfill \null \\[0.5em]
\end{lem}

\begin{pf}
Assume that for all $k>0$, there exists $x_k$ in $\C$ and $s_k$ in $\sign(f(x_k))\subset[-1,1]$ such that
\\[0.5em]\null \hfill $\displaystyle 
0>\eta(x_k,s_k) \geq -\frac{1}{k}  \ .
$\hfill \null \\[0.5em]
Then, $\eta(x_k,s_k)$ tends to $0$ when $k$ tends to infinity. Besides, there exists a subsequence $(k_m)$ such that $x_{k_m}$ tends to $x^*$ in $\C$ and $s_{k_m}$ tends to $s^*$ in $[-1,1]$. Since $\eta$ is continuous, it follows that $\eta(x^*,s^*)=0$ and we will have a contradiction if $s^*\in\sign(f(x^*))$. If $f(x^*)$ is not zero, then by continuity of $f$, $s^*$ is equal to the sign of  $f(x^*)$, and otherwise, $s^*\in[-1,1]=\sign(f(x^*))$. Thus, $s^*\in\sign(f(x^*))$ in all cases.
\end{pf}

\begin{lem}
\label{lem2}
Let $\eta$ be a function defined on $\RR^n$ homogeneous with degree $d$ and weight vector $r=(r_1,...,r_n)$, and $V$ a positive definite proper function defined on $\RR^n$ homogeneous of degree $d_V$ with same weight vector $r$. Define $\C=V^{-1}(\{1\})$.  If there exists $\alpha$ such that for all $x$ in $\C$
$$
\eta(x) < \alpha \ ,
$$
then for all $x$ in $\RR^n\setminus\{0\}$,
$\displaystyle
\eta(x) < \alpha V(x)^{\frac{d}{d_V}}
$.
\end{lem}

\begin{pf}
Let $x$ in $\RR^n\setminus\{0\}$. We have
$\displaystyle
\bar x = \frac{x_i}{V(x)^{\frac{r_i}{d_V}}} 
$ in $\C$. 
Thus
$
\eta(\bar x)  <\alpha
$
and by homogeneity
\\[0.5em]\null \hfill $\displaystyle
\frac{1}{V(x)^{\frac{d}{d_V}}}\eta( x)  < \alpha
$\hfill \null \\[0.5em]
which gives the required inequality.
\end{pf}

\begin{lem}
\label{lem}
Let $\eta$ be a homogeneous function of degree $d$ and weight vector $r$ defined on $\RR^n$ by
\\[0.5em]\null \hfill $\displaystyle 
\eta(x)=\max_{s\in\sign(f(x))} \tilde{\eta}(x,s) 
$\hfill \null \\[0.5em]
where $\tilde{\eta}$ is a continuous function defined on $\RR^{n+1}$ and $f$ a
continuous function defined on $\RR^n$. Consider a continuous function $\gamma$
homogeneous with same degree and weight vector such that, for all $x$ in $\RR^n\setminus\{0\}$
and $s$ in $\sign(f(x))$
\\[0.5em]\null \hfill $
\begin{array}{lcl}
\gamma(x) \geq 0 \ , \\
\gamma(x) = 0 \quad \Rightarrow \quad \tilde\eta(x,s) < 0 \ .
\end{array}
$\hfill \null \\[0.5em]
Then, there exists $k_0>0$ such that, for all $x$ in $\RR^n\setminus\{0\}$,
\\[0.5em]\null \hfill $
\eta(x) - k_0 \ \gamma(x) <0 \ .
$\hfill \null 
\end{lem}

\begin{pf}
	Define the homogeneous definite positive function 
	$\displaystyle
	V(x) = \sum_{i=1}^n{|x_i|^{\frac{d}{r_i}}}
	$
	and consider the compact set 
	$
	\C = V^{-1}(\{1\})
	$.
	Assume that for all $k>0$, there exists $x_k$ in $\C$ and $s_k$ in $\sign(f(x_k))$ such that
	$$
	\tilde\eta(x_k,s_k) \geq k \ \gamma(x_k) \geq 0
	$$
	$\tilde\eta$ is continuous, and thus bounded on the compact set $\C\times[-1,1]$. Therefore, $\gamma(x_k)$ tends to $0$ when $k$ tends to infinity. Besides, there exists a subsequence $(k_m)$ such that $x_{k_m}$ tends to $x^*$ in $\C$ and $s_{k_m}$ tends to $s^*$ in $[-1,1]$. It follows that $\gamma(x^*)=0$ since $\gamma$ is continuous. But with the same argument as in the proof of Lemma \ref{lemBis}, we have $s^*\in\sign(f(x^*))$. It yields that $\tilde\eta(x^*,s^*)<0$ by assumption and we have a contradiction.
	
	Therefore, there exists $k_0$ such that 
	\\[0.6em]\null \hfill $\displaystyle 
	\tilde\eta(x,s) - k_0 \ \gamma(x) <0
	$\hfill \null \\[0.6em]
	for all $x$ in $\C$ and all $s$ in $\sign(f(x))$. Thus, with Lemma \ref{lemBis} there exists $\alpha>0$ such that  
	\\[0.7em]\null \hfill $\displaystyle 
	\tilde\eta(x,s) - k_0 \ \gamma(x) \leq-\alpha
	$\hfill \null \\[0.3em]
	so that
	$$
	\eta(x) - k_0 \ \gamma(x) <0
	$$
	for any $x$ in $\C$. The result follows applying Lemma \ref{lem2}.
\end{pf}

\startrevision
\begin{lem}\label{lem_LyapIss}
For a positive bounded continuous function $t\mapsto  c(t)$ and an absolutely continuous
function $t\mapsto  \nu (t)$ satisfying
\\[0.3em]\null \hfill for almost all $t$ such that $\nu (t) \geq c(t)$ then $\dot \nu (t) \leq -\nu (t)^d$
\hfill \null \\[0.3em]
with $d$ in $]0,1[$.
Then, for all $t$ in $[0,T[$
\\[0.5em]$\displaystyle 
\nu (t) \leq  \max\left\{0,\max\{\nu (0)-c(0),0\}^{1-d}-t\right\}^{1/(1-d)} 
$\hfill \null \\\null \hfill  $\displaystyle 
+ \sup_{s\in [0,t]}c(s)\  .$
\end{lem}
\begin{pf}
This is a direct consequence of the fact that we have for almost all $t$ such that $\max\{\nu (t)-c,0\}) $ is $C^1$
\\[0.5em]\null \hfill $\displaystyle 
\dot{\overparen{\max\{\nu (t)-c,0\}}} \leq - \max\{\nu (t)-c,0\}^d
$\hfill \null 
\end{pf}
\stoprevision

\bibliography{biblio}      
\bibliographystyle{plain}

\end{document}